\documentclass[11pt]{article}
\usepackage{amsmath, amsthm, amssymb, eucal}
\newcommand{\partentry}[1]{\addtocontents{toc}
{\small\bfseries#1\hfill\thepage\par}}
\makeatletter
\def\@part[#1]#2{%
    \ifnum \c@secnumdepth >\m@ne
      \refstepcounter{part}
      \partentry{\protect\makebox[2em][l]{\thepart}#1}
\else
      \partentry{#1}
    \fi
    {\noindent\normalfont\Large\bfseries\thepart\hspace{1em}#2\par}
    \nobreak
    \vskip 3ex
    \@afterheading}
\def\@spart#1{%
    {\noindent\normalfont\Large\bfseries #1\par} 
     \nobreak
     \vskip 3ex
     \@afterheading}
\renewcommand\section{\@startsection{section}{1}{\z@}
{-3.5ex \@plus -1ex \@minus -.2ex}
{2ex \@plus.2ex}
{\large\bfseries}}
\renewcommand\subsection{
\@ifstar{\setcounter{subsection}{\value{equation}}
\@startsection{subsection}{2}{\z@}
	{1.75ex \@plus.5ex \@minus.2ex}
 	{-.4em} 
	{\itshape}*}
{\setcounter{subsection}{\value{equation}}
\stepcounter{equation}
\@startsection{subsection}{2}{\z@}
	{1.75ex \@plus.5ex \@minus.2ex}
 	{-.4em} 
	{\itshape}}}
\def\@seccntformat#1{\@ifundefined{#1@cntformat}
{\csname the#1\endcsname\quad}
{\csname #1@cntformat\endcsname}}
\def\section@cntformat{\thesection.~}
\def\subsection@cntformat{(\thesubsection)\ }
\renewcommand*\l@section{\mdseries\small\@dottedtocline{1}{1.5em}{2em}}
\makeatother

\numberwithin{equation}{section}
\theoremstyle{plain}
\newcounter{mainth}

\newtheorem{maintheorem}[mainth]{Theorem}
\swapnumbers
\newtheorem{theorem}[equation]{Theorem}

\newtheorem{lemma}[equation]{Lemma}
\newtheorem{proposition}[equation]{Proposition}
\theoremstyle{definition}
\newtheorem{definition}[equation]{Definition}
\theoremstyle{remark}
\newtheorem{remark}[equation]{Remark}

\newcommand{\cG}{\mathcal{G}}
\newcommand{\cH}{\mathcal{H}}

\newcommand{\cO}{\mathcal{O}}
\newcommand{\cP}{\mathcal{P}}

\newcommand{\cS}{\mathcal{S}}
\newcommand{\cV}{\mathcal{V}}
\newcommand{\frb}{\mathfrak{b}}
\newcommand{\frg}{\mathfrak{g}}
\newcommand{\frn}{\mathfrak{n}}
\newcommand{\frh}{\mathfrak{h}}

\newcommand{\frB}{\mathfrak{B}}
\newcommand{\frL}{\mathfrak{L}}
\newcommand{\frM}{\mathfrak{M}}
\newcommand{\bA}{\mathbb{A}}
\newcommand{\bC}{\mathbb{C}}
\newcommand{\bH}{\mathbb{H}}
\newcommand{\bP}{\mathbb{P}}
\newcommand{\bQ}{\mathbb{Q}}

\newcommand{\bZ}{\mathbb{Z}}
\newcommand{\bfH}{\mathbf{H}}
\newcommand{\bfX}{\mathbf{X}}
\newcommand{\bfU}{\mathbf{U}}
\newcommand{\me}{\mathrm{e}}
\newcommand{\vep}{\varepsilon}
\newcommand{\SL}{\mathrm{SL}}

\newcommand{\ad}{\mathrm{ad}}
\newcommand{\Ext}{\mathrm{Ext}}
\newcommand{\gen}{\mathrm{Gen}}
\newcommand{\Hom}{\mathrm{Hom}}
\newcommand{\prim}{\mathrm{Prim}\,}
\newcommand{\sym}{\mathrm{S}\mspace{1mu}}
\newcommand{\Tor}{\mathrm{Tor}}
\newcommand{\Tr}{\mathrm{Tr}}
\newcommand{\waff}{W_\mathrm{aff}}
\newcommand{\wD}{\widehat{D}}
\newcommand{\im}{\mathrm{Im}\,}
\newcommand{\Gr}{\mathrm{Gr}}
\newcommand{\oSig}{{\overline{\Sigma}}}
\newcommand{\ofrac}[2]{\genfrac{}{}{0pt}{}{#1}{#2}}

\begin{document}

\title{\textbf{The strong Macdonald conjecture and\\
		Hodge theory on the Loop Grassmannian}}
\author{Susanna Fishel \and Ian Grojnowski \and Constantin Teleman} 
\date{\today}
\maketitle

\begin{quote}
\abstract{\noindent We prove the strong Macdonald conjecture of Hanlon 
and Feigin for reductive groups $G$. In a geometric reformulation, we 
show that the Dolbeault cohomology $H^q(X;\Omega^p)$ of the loop 
Grassmannian $X$ is freely generated by de Rham's forms on the disk 
coupled to algebra generators of $H^\bullet (BG)$. Equating Euler 
characteristics of the two gives an identity, independently known to 
Macdonald \cite{mac}, which generalises Ramanujan's ${}_1\psi_1$ sum. 
For simply laced root systems at level $1$, we find a `strong form' 
of Bailey's ${}_4\psi_4$ sum. Failure of Hodge decomposition implies 
the \textit{singularity} of $X$, and of the algebraic loop groups. 
Some of our results were announced in \cite{tel2}.}
\end{quote}

\section*{Introduction}
This article address some basic questions concerning the cohomology of 
affine Lie algebras and their flag varieties. Its chapters are closely 
related, but have somewhat different flavours, and the methods used may 
well appeal to different readers. Chapter I proves the \textit{strong 
Macdonald constant term conjectures} of Hanlon \cite{han1} and Feigin 
\cite{feig1}, describing the cohomologies of the Lie algebras $\frg[z]/z^n$ 
of truncated polynomials with values in a reductive Lie algebra $\frg$ 
and of the \textit{graded} Lie algebra $\frg[z,s]$ of $\frg$-valued skew 
polynomials in an even variable $z$ and an odd one $s$ (Theorems \ref
{trunc} and \ref{sym}). The proof uses little more than linear algebra, 
and, while Nakano's identity (\ref{315}) effects a substantial 
simplification, we have included a brutal computational by-pass 
in Appendix \ref{app}, to avoid reliance on external sources. 

Chapter II discusses the Dolbeault cohomology $H^q(\Omega^p)$ of flag 
varieties of loop groups. In addition to the ``Macdonald cohomology", the 
methods and proofs draw heavily on \cite{tel3}. For the \textit{loop 
Grassmannian} $X:= G((z))/G[[z]]$, we obtain the free algebra generated 
by copies of the spaces $\bC[[z]]$ and $\bC[[z]]dz$, in bi-degrees $(p,q) 
=(m,m)$, respectively $(m,m+1)$, as $m$ ranges over the exponents of $\frg$. 
Moreover, de Rham's operator $\partial: H^q (\Omega^p)\to H^q(\Omega^{p+1})$ 
is induced by $d: \bC[[z]] \to\bC[[z]]dz$ on generators.  

A noteworthy consequence of our computation is the \textit{failure of Hodge 
decomposition}, 
\[
H^n(X;\bC)\neq\bigoplus\nolimits_{p+q=n} H^q(X;\Omega^p). 
\]
Because $X$ is a union of projective varieties, this implies that $X$ 
\textit{is not smooth}, in the sense that it is not locally expressible as 
an increasing union of smooth complex-analytic sub-varieties (Theorem \ref
{65}). We are thus dealing with a \textit{homogeneous variety} which is 
\textit{singular everywhere}. We are unable to offer a true geometric 
explanation of this striking fact. 

Our results generalise to an arbitrary smooth affine curve $\Sigma$. The 
Macdonald cohomology involves now the Lie algebra $\frg[\Sigma,s]$ of 
$\frg[s]$-valued algebraic maps, while $X$ is replaced by the \textit
{thick flag variety} $\bfX_\Sigma$ of \S\ref{7}. In this generality, the 
question requires more insight than is provided by the listing of co-cycles 
in Theorem \ref{sym}. Thus, after re-interpreting the Macdonald cohomology 
as the (algebraic) Dolbeault cohomology of the classifying stack $BG[[z]]$, 
and the flag varieties $\bfX_\Sigma$ as moduli of $G$-bundles on $\Sigma$ 
trivialised near $\infty$, \S\ref{unified} gives a uniform construction of 
all generating Dolbeault classes. Inspired by the Atiyah-Bott description 
of the cohomology generators for the moduli of $G$-bundles, our construction 
is a Dolbeault refinement, based on the Atiyah class of the universal bundle 
and invariant polynomials on $\frg$, in lieu of the Chern classes.

The more geometric perspective leads us to study $H^q(X; \Omega^p \otimes\cV)$, 
for certain vector bundles $\cV$; this ushers in Chapter III. In \S\ref
{twist}, we find a beautiful answer for simply laced groups and the 
level $1$ line bundle $\cO(1)$. In general, we can define, for each level 
$h\ge 0$ and $G$-representation $V$, the formal Euler series in $t$ and $z$ 
with coefficients in the character ring of $G$:
\[
P_{h,V} = \sum\nolimits_{p,q} (-1)^q (-t)^p \mathrm{ch}\, 
		H^q\left(X;\Omega^p(h) \otimes\cV\right),
\] 
where the vector bundle $\cV$ is associated to the $G$-module $V$ as in 
\S\ref{11-7} and $z$ carries the weights of the $\bC^\times$-scaling action 
on $X$. These series, expressible using the Kac character formula, are 
affine analogues of the Hall-Littlewood symmetric functions, and their 
complexity leaves little hope for an explicit description of the cohomologies. 
On the other hand, the finite Hall-Littlewood functions are related to 
certain filtrations on weight spaces of $G$-modules, studied by Kostant, 
Lusztig and Ran\'ee Brylinski in general. We find in \S\ref{afhall} that 
such a relationship persists in the affine case \textit{at positive level}. 
Failure of the level zero theory is captured precisely by the Macdonald 
cohomology, or by its Dolbeault counterpart in Chapter II; whereas the 
good behaviour at positive level relies on a higher-cohomology vanishing 
(Theorem~\ref{vanish}). 

We emphasise that finite-dimensional analogues of our results (Remarks 
\ref{11-0} and \ref{11-10}), which are known to carry geometric information 
about the $G$-flag variety $G/B$ and the nilpotent cone in $\frg$, can be 
deduced from standard Hodge theory or other cohomology vanishing results 
(such as the Grauert-Riemenschneider theorem, applied to the moment map 
$\mu: T^*G/B\to\frg^*$). No such general theorems are available in the loop 
group case; our results provide a substitute for this. Developing the full 
theory would take us too far afield, and we postpone it to a future paper, 
but \S\ref{11} illustrates it with a simple example. 

Finally, just as the strong Macdonald conjecture refines a combinatorial 
identity, our new results also have combinatorial applications. Comparing 
our answer for $H^q(X;\Omega^p(h))$ with the Kac character formula for 
$P_{h,\bC}$ leads to $q$-basic hyper-geometric summation identities. For 
$\SL_2$, this is a specialisation of Ramanujan's ${}_1\psi_1$ sum. For 
general affine root systems, these identities were independently discovered 
by Macdonald \cite{mac}. The level one identity for $\SL_2$ comes from a 
specialised Bailey ${}_4\psi_4$ sum; its extension to simply laced root 
systems seems new. 

Most of the work for this paper dates back to 1998, and the authors have 
lectured on it at various times; the original announcement is in \cite
{tel2}, and a more leisurely survey is \cite{gro}. We apologise for the 
delay in preparing the final version.

\vskip4ex
\noindent\small\textit{Acknowledgements.} The first substantial portion 
of this paper (Chapter I) was written and circulated in 2001, during the 
most enjoyable programme on ``Symmetric Functions and Macdonald Polynomials" 
at the Newton Institute in Cambridge. We wish to thank numerous colleagues, 
among whom are E.~Frenkel, P.~Hanlon, S.~Kumar, I.G.~Macdonald, S.~Milne, 
for their comments and interest, as well as their patience. The third author 
was originally supported by an NSF Postdoctoral Fellowship. 
\normalsize

\vskip4ex
\begin{minipage}[t]{12cm}
\tableofcontents
\end{minipage}
\vskip4ex
\section*{Definitions and Notation.} 
Our (Lie) algebras and vector spaces are defined over $\bC$. Certain 
vector spaces, such as $\bC[[z]]$, have natural inverse limit topologies, 
and $*$ will then denote their \textit{continuous} duals; this way, 
$\bC[[z]]^{**}\cong\bC[[z]]$. Completed tensor products or powers of such 
spaces will be indicated by $\widehat\otimes$, $\hat\sym{}^p$, $\hat 
\Lambda{}^p$.
 
\subsection{Lie algebra (co)homology.}\label{26}
The \textit{Lie algebra homology Koszul complex} \cite{kos} of a Lie 
algebra $\frL$ with coefficients in a module $V$ is $\Lambda^\bullet 
\frL\otimes V$, homologically graded by $\bullet$, with differential
\[\begin{split}
\delta (\lambda_1\wedge\ldots\wedge\lambda_n\otimes v)
&=\sum\nolimits_p (-1)^p\lambda_1\wedge\ldots\wedge\hat\lambda_p\wedge
\ldots \wedge\lambda_n\otimes \lambda_p(v)\\
&+{}\sum\nolimits_{p<q}(-1)^{p+q}[\lambda_p,\lambda_q]\wedge\lambda_1
\wedge\ldots\wedge\hat\lambda_p\wedge\ldots\wedge \hat\lambda_q\wedge
\ldots\wedge\lambda_n\otimes v;
\end{split}\]
hats indicate missing factors. Its homology $H_\bullet(\frL;V)$ is the 
\textit{Lie algebra homology of $\frL$ with coefficients in} $V$. If 
$\frg\subseteq\frL$ is a sub-algebra, $\delta$ descends to the quotient 
$\left(\Lambda(\frL/\frg)\otimes V\right)/\frg\left(\Lambda(\frL/\frg)
\otimes V\right)$ of co-invariants under $\frg$, which resolves the \textit
{relative homology} $H_\bullet(\frL,\frg;V)$. We denote by $H_\bullet(\frL)$ 
the homology with coefficients in the trivial one-dimensional module. 

Dual to these are the \textit{cohomology complexes}, with underlying 
spaces $\Hom(\Lambda^\bullet\frL;W)$; the cohomology is denoted 
$H^\bullet(\frL;W)$, or $H^\bullet(\frL,\frg;W)$ in the relative case.
They are the full duals of the homologies, when $W$ is the full dual of 
$V$. If $W$ is an algebra and $\frL$ acts by derivations, the Koszul complex 
is a differential graded algebra. Similarly, the homology complex is a 
differential graded co-algebra, when $V$ is a co-algebra and $\frL$ acts 
by co-derivations.

\begin{remark}
More abstractly, $H_k(\frL;V)=\Tor^\frL_k(\bC;V)$ and $H^k(\frL;V) =
\Ext_\frL^k(\bC;V)$ in the category of $\frL$-modules. If $\frg\subseteq
\frL$ is reductive, and $\frL$ (via ad) and $V$ are semi-simple $\frg
$-modules, the relative homologies are the $\Tor$ groups in the category 
of \textit{$\frg$-semi-simple} $\frL$-modules.
\end{remark}

\subsection{Exponents.}\label{exp}
Either of the following statements defines the \textit{exponents} $m_1,
\ldots ,m_\ell$ of a reductive Lie algebra $\frg$ of rank $\ell$:
\begin{itemize}
\item the algebra $(\sym\frg^*)^\frg$ of polynomials on $\frg$ which are
invariant under the co-adjoint action is a free symmetric algebra generated 
in degrees $m_1+1,\ldots ,m_\ell +1$;
\item the sub-algebra $(\Lambda\frg)^\frg$ of ad-invariants in the full 
exterior algebra of $\frg$ is a free exterior algebra generated in degrees
$2m_1+1,\ldots,2m_\ell +1$. 
\end{itemize}
For instance, when $\frg =\mathfrak{gl}_n$, $\ell=n$ and $\left(m_1, 
\ldots,m_n\right) =\left(0,\ldots,n-1\right)$. The first algebra is 
also naturally isomorphic to the cohomology $H^\bullet(BG;\bC)$, if we 
set $\deg\frg =2$.

\subsection{Generators.} Most cohomologies in this paper will be free 
graded polynomial (or power series) algebras, which are canonically 
described by identifying their spaces of indecomposables\footnote{Recall 
that the space of indecomposables of a non-negatively graded algebra 
$A^{\bullet}$ is $A^{>0}/(A^{>0})^2$. If $A^\bullet$ is a free algebra 
over $A^0$, a graded $A^0$-lifting of the indecomposables in $A^\bullet$ 
gives a space of algebra generators.} with those for $H^\bullet(BG)$, 
tensored with suitable graded vector spaces $V^\bullet$ (cf.\ Theorem \ref{sym}). However, we can choose once and for all a space $\gen^
\bullet(BG)$ spanned by homogeneous free generators for the cohomology, 
and identify our cohomologies as the free algebras on $\gen^\bullet(BG)
\otimes V^\bullet$. There are many choices of generators,\footnote{Natural 
examples for $GL_n$ include the Chern classes and the traces $\Tr F^k$ 
of the universal curvature form $F$.} but our explicit constructions of 
cohomology classes from invariant polynomials serves to make this 
second description canonical.

\subsection{Fourier basis.} \label{02} When $G$ is semi-simple, we will 
choose a compact form and a basis of self-adjoint elements $\xi_a$ in 
$\frg$, orthonormal in the Killing form. Call, for $m\ge 0$, $\psi^a(-m)$ 
and $\sigma^a(-m)$ the elements of $\Lambda^1\frg[z]^*$ and $\sym^1\frg[z]
^*$ dual to the basis $z^m\cdot\xi_a$ of the Lie algebra $\frg[z]$. We 
abusively write $\xi_{[a,b]}$ for $[\xi_a, \xi_b]$, and similarly 
$\psi^{[a,b]}(m)$ for $\ad^*_{\xi_a}\psi^b(m)$, etc.

\part{The strong Macdonald conjecture} 
\section{Statements}\label{1} 
\subsection{Background.}
The strong Macdonald conjectures describe the cohomologies of the truncated 
Lie algebras $\frg[z]/z^n$ and of the graded Lie algebra $\frg[z,s]$. The 
first conjecture is due to Hanlon \cite{han1}, who also proved it for 
$\mathfrak{gl}_n$ \cite{han2}. The conjecture may have been independently 
known to Feigin \cite{feig1}, who in \cite{feig2} related it to the 
cohomology of $\frg[z,s]$. Feigin also outlined a computation of the latter; 
but we are unsure whether it can be carried out as indicated.\footnote{One 
particular step, the lemma on p.~93 of \cite{feig2}, seems incorrect: the 
analogous statement fails for absolute cohomology when $Q=\partial/\partial
\xi$, and nothing in the suggested argument seems to account for that.} 
While we could not fill the gap, we do confirm the conjectures by a 
different route: we compute the cohomology of $\frg[z,s]$ by finding the 
\textit{harmonic co-cycles} in the Koszul complex, in a suitable metric. 
Feigin's argument then recovers the cohomology of the truncated Lie algebra.

The success of our Laplacian approach relies on the specific metric used
on the Koszul complex and originates in the K\"ahler geometry of the loop 
Grassmannian. The latter is responsible for an identity between two 
different Laplacians, far from obvious in Lie algebra form, which implies 
here that the harmonic co-cycles form a \textit{sub-algebra} and allows 
their computation. We do not know of a computation using the more obvious 
Killing metric; its harmonic co-cycles are not closed under multiplication. 

\subsection*{Truncated algebras.} The following affirms Hanlon's original 
conjecture for reductive $\frg$. Note that the cohomology of $\frg[z]/z^n$ 
decomposes by $z$-weight, in addition to the ordinary grading.
\begin{maintheorem}\label{trunc}  
$H^\bullet(\frg[z]/z^n)$ is a free exterior algebra on $n\cdot\ell$ 
generators, with $n$ generators in cohomology degree $2m+1$ and $z$-weights 
equal to the negatives of $0,mn+1,mn+2,\ldots,mn+n-1$, for each exponent 
$m = m_1,\ldots,m_\ell$. 
\end{maintheorem}

\begin{remark}\label{12}\begin{trivlist}\itemsep0.2ex 
\item(i) Ignoring  $z$-weights leads to an \textit{abstract} ring 
isomorphism $H^\bullet(\frg[z]/z^n)\cong H^\bullet(\frg)^{\otimes n}$.
\item(ii) The degree-wise lower bound $\dim H^\bullet(\frg[z]/z^n)\ge 
\dim H^\bullet(\frg)^{\otimes n}$ holds for \textit{any} Lie algebra 
$\frg$. Namely, $\frg[z]/z^n$ is a degeneration of $\frg[z]/(z^n-\vep)$, 
as $\vep\to 0$. When $\vep\ne 0$, the quotient is isomorphic to $\frg^
{\oplus n}$, whose cohomology is $H^\bullet(\frg)^{\otimes n}$, and the 
ranks are upper semi-continuous. However, this argument says nothing about 
the \textit{ring} structure. 
\item(iii) There is a \textit{natural} factorisation $H^\bullet(\frg[z]/
z^n) = H^\bullet(\frg)\otimes H^\bullet(\frg[z]/z^n,\frg)$, and the first 
factor has $z$-weight $0$. Indeed, reductivity of $\frg$ leads to a spectral 
sequence \cite{kos} with
\[
E_2^{p,q} = H^q(\frg)\otimes H^p(\frg[z]/z^n,\frg) \Rightarrow 
		H^{p+q}(\frg[z]/z^n),
\]
whose collapse there is secured by the evaluation map $\frg[z]/z^n \to 
\frg$, which provides a lifting of the left edge $H^q(\frg)$ in the 
abutment and denies the possibility of higher differentials. 
\end{trivlist}
\end{remark}

\subsection{Relation to cyclic homology.} A conceptual formulation of 
Theorem \ref{trunc} was suggested independently by Feigin and Loday. Given 
a skew-commutative algebra  $A$ and any Lie algebra $\frg$, an invariant 
polynomial $\Phi$ of degree $(m+1)$ on $\frg$ determines a linear map from 
the dual of $HC_n^{(m)}(A)$, the $m$th \textit{Adams component} of the $n$th 
\textit{cyclic homology} group of $A$, to $H^{n+1}(\frg\otimes A)$ (see our 
Theorem \ref{sym} for the case of interest here, and \cite{tel2} (2.2), or 
the comprehensive discussion in \cite{lod} in general). When  $\frg$ is 
reductive, Loday suggested that these maps might be injective, and that 
$H^\bullet(\frg\otimes A)$ might be freely generated by their images, as 
$\Phi$ ranges over a set of generators of the ring of invariant polynomials. 
The Adams degree $m$ will then range over the exponents $m_1,\ldots,m_\ell$. 
Thus, for $A=\bC$, $HC_n^{(m)} = 0$ for $n\ne 2m$, while $HC_{2m}^{(m)} = 
\bC$; and we recover the well-known description of $H^\bullet(\frg)$. For 
$\frg=\mathfrak{gl}_\infty$ and any associative, unital, graded $A$, this 
is the theorem of Loday-Quillen \cite{loquil} and Tsygan \cite{ts}. It 
emerges from its proof that Theorem \ref{trunc} affirms Loday's conjecture 
for $\bC[z]/z^n$, while \eqref{14} below does the same for the graded 
algebra $\bC[z,s]$. (The conjecture fails in general \cite{tel2}.) 

\subsection{The super-algebra.} 
The graded space $\frg[z,s]$ of $\frg$-valued skew polynomials in $z$ and 
$s$, with $\deg z=0$ and $\deg s = 1$, is an infinite-dimensional \textit 
{graded} Lie algebra, isomorphic to the semi-direct product $\frg[z]\ltimes 
s\frg[z]$ (for the adjoint action), with zero bracket in the second factor. 
We shall give three increasingly concrete descriptions \eqref{14}, \eqref
{19}, \eqref{sym} for its (co)homology. We start with \textit{homology}, 
which has a natural \textit{co-algebra} structure. As in Remark \ref{12}.iii, 
we factor $H_\bullet(\frg[z,s])$ as $H_\bullet(\frg)\otimes H_\bullet 
(\frg[z,s],\frg)$; the first factor behaves rather differently from the 
rest, and is best set aside. 

\begin{theorem}\label{14} 
$H_\bullet(\frg[z,s],\frg)$ is isomorphic to the free, graded co-commutative 
co-algebra whose space of primitives is the direct sum of copies of 
$\bC[z]\cdot s^{\otimes(m+1)}$, in total degree $2m+2$, and of $\bC[z]dz
\cdot s^{\otimes m}$, in total degree $2m+1$, as m ranges over the exponents 
$m_1,\ldots,m_\ell$. The isomorphism respects $(z,s)$-weights. 
\end{theorem}

\begin{remark}\label{15}
\begin{trivlist}\itemsep0ex 
\item(i) The \textit{total degree} $\bullet$ includes that of $s$. As 
multi-linear tensors in $\frg[z,s]$, both types of cycles have degree $m+1$. 
\item(ii) A free co-commutative co-algebra is isomorphic, as a vector space,  
to the graded symmetric algebra on its primitives; but there is no 
\textit{a priori} algebra structure on homology. 
\end{trivlist}\end{remark}

The description \eqref{14} is not quite canonical. If $P_{(k)}$ is the 
space of $k$th degree primitives in the quotient co-algebra $\sym\frg/
[\frg,\sym\frg]$, canonical descriptions of our primitives 
are \begin{equation}\begin{split}\label{16}
\bigoplus\nolimits_m &P_{(m+1)}\otimes \bC[z]\cdot s(ds)^m, \\
\bigoplus\nolimits_m &P_{(m+1)}\otimes
	\frac{\bC[z]\cdot(ds)^m+\bC[z]dz\cdot s(ds)^{m-1}}{
d\left(\bC[z]\cdot s(ds)^{m-1}\right)}.
\end{split}\end{equation}
The right factors are the cyclic homology components $HC_{2m+1}^{(m)}$
and $HC_{2m}^{(m)}$ of the \textit{non-unital} algebra $\bC[z,s]\ominus\bC$.
The last factor, $HC_{2m}^{(m)}$, is identifiable with $\bC[z]dz\cdot
s(ds)^{m-1}$, for $m\ne 0$, and with $\bC[z]/\bC$ if $m=0$. This description 
is compatible with the action of \textit{super-vector fields} in $z$ and $s$ 
(see Remark \ref{255} below), whereas \eqref{14} only captures the action 
of vector fields in $z$.

\subsection{Restatement without super-algebras.}\label{24} 
There is a natural isomorphism between $H_\bullet(\frL; \Lambda^\bullet V)$ 
and the homology of the semi-direct product Lie algebra $\frL\ltimes V$, 
with zero bracket on $V$ \cite{kos}. Its graded version, applied to $\frL 
= \frg[z]$ and the \textit{odd} vector space $V = s\frg[z]$, is the 
equality 
\begin{equation}\label{18}			
H_n(\frg[z,s],\frg) = \bigoplus_{p+q=n}H_{q-p}\left(\frg[z],\frg;
\sym^p(s\frg[z])\right); 
\end{equation}
note that elements of $s\frg[z]$ carry homology degree $2$ (Remark \ref
{15}.i). We can restate Theorem \ref{14} as follows:

\begin{theorem}\label{19} 
$H_\bullet\left(\frg[z],\frg;\sym(s\frg[z])\right)$ is isomorphic to the 
free graded co-commutative co-algebra with primitive space $\bC[z]\cdot 
s^{\otimes (m+1)}$, in degree $0$, and primitive space $\bC[z]dz\cdot 
s^{\otimes m}$ in degree $1$, as $m$ ranges over the exponents $m_1,\ldots, 
m_\ell$. The isomorphism preserves $z$-and $s$-weights. 
\end{theorem}

\subsection{Cohomology.} \label{1-11}
While $H^\bullet(\frg[z,s],\frg)$ is obtained from \eqref{18} by duality, 
infinite-dimensionality makes it a bit awkward, and we opt for a 
\textit{restricted duality}, defined using the direct sum of the 
$(s,z)$-weight spaces in the dual of the Koszul complex \eqref{26}. These 
weight spaces are finite-dimensional and are preserved by the Koszul 
differential. The resulting restricted Lie algebra cohomology $H_{res}^
\bullet(\frg[z,s],\frg)$ is the direct sum of weight spaces in the full 
dual of \eqref{18}.

\begin{maintheorem} \label{sym} 
$H_{res}^\bullet\left(\frg[z],\frg;\sym\frg[z]^*\right)$ is isomorphic to 
the free graded commutative algebra generated by the restricted duals of 
$\bigoplus_mP_{(m+1)}\otimes \bC[z]$ and $\bigoplus_mP_{(m+1)}\otimes 
\bC[z]dz$, in cohomology degrees $0$ and $1$ and symmetric degrees $m+1$ 
and $m$, respectively. \\
Specifically, an invariant linear map $\Phi: \sym^{m+1}\frg \to \bC$ 
determines linear maps 
\begin{align*}
S_\Phi:\sym^{m+1}\frg[z]&\to\bC[z], \\	
	&\sigma_0\cdot\sigma_1\cdot\ldots\cdot\sigma_m\mapsto
		\Phi\left(\sigma_0(z),\sigma_1(z),\ldots,\sigma_m(z)\right)\\
E_\Phi:\Lambda^1\left(\frg[z]/\frg\right)&\otimes\sym^m\frg[z]\to \bC[z]dz, \\
	&\psi\otimes \sigma_1\cdot \ldots \cdot \sigma_m\mapsto
		\Phi\left(d\psi(z),\sigma_1(z),\ldots,\sigma_m(z)\right).
\end{align*}
The coefficients $S_\Phi(-n)$, $E_\Phi(-n)$ of $z^n$, resp. $z^{n-1}dz$ are 
restricted $0$- and $1$-cocycles and $H_{res}^\bullet$ is freely generated 
by these, as $\Phi$ ranges over a generating set of invariant polynomials 
on $\frg$.
\end{maintheorem}

\noindent To illustrate, here are the cocycles associated to the Killing form
on $\frg$ (notations as in \S\ref{02}):
\[
S(-n) = \sum_{\ofrac{1\le a\le\dim G}{0\le p\le n}}
			\sigma^a(-p)\sigma^a(p-n), \qquad	
E(-n) = \sum_{\ofrac{1\le a\le\dim G}{0<p\le n}}
			p\psi^a(-p)\sigma^a(p-n).
\]

We close this section with two generalisations of Theorem \ref{sym}. The 
first will be proved in \S\ref{4}; the second relies on more difficult 
techniques, and will only be proved in \S\ref{10}.

\subsection{The Iwahori sub-algebra.} 
Here, we replace $\frg[z]$ with an \textit{Iwahori sub-algebra} $\frB\subset 
\frg[z]$, the inverse image of a Borel sub-algebra $\frb\subset\frg$ under 
the evaluation at $z=0$. The cocycles $S_\Phi(0)$ generate a copy of 
$(\sym^\bullet \frg^*)^\frg$ within $H_{res}^{2\bullet}\left( \frg[z,s], 
\frg\right)$. With $\frh:= \frb/[\frb,\frb]$, isomorphic to a Cartan 
sub-algebra, a similar inclusion $\sym^\bullet\frh^*\to H_{res}^ {2\bullet} 
\left(\frB[s],\frh\right)$ results from identifying $\frh^*$ with the $\frB
$-invariants in $\frB^*$. Recall that $(\sym\frg^*)^\frg$ embeds in $\sym 
\frh^*$ (as the Weyl-invariant sub-algebra). It turns out that, when passing 
from $\frg[z]$ to $\frB$, the factor $(\sym\frg^*)^\frg$ is replaced with 
$\sym\frh^*$.

\begin{theorem}\label{111} 
$H_{res}^\bullet\left(\frB[s],\frh\right) \cong
H_{res}^\bullet\left(\frg[z,s],\frg\right)\otimes_{(\sym(s\frg)^*)^G}
\sym(s\frh)^*$.
\end{theorem}

\subsection{Affine curves.} \label{arbaff}
Our second generalisation replaces $\frg[z]$ by the $\frg$-valued algebraic 
functions on a smooth affine curve $\Sigma$. The space $\frg[\Sigma]$ has 
no restricted dual as in \S\ref{1-11}, so we use full duals in the Koszul 
complex; consequently, the cohomology will be a power series algebra. 
Moreover, there is now a contribution from the cohomology with  constant 
coefficients (whereas before we had $H^\bullet(\frg[z],\frg;\bC)=\bC$, by 
\cite{garlep}). The last cohomology is described in \ref{10-6}.

\begin{theorem}\label{1-15}
For a smooth affine curve $\Sigma$, the cohomology $H^\bullet\left(\frg
[\Sigma];(\sym\frg[\Sigma])^*\right)$ is densely generated over 
$H^\bullet(\frg[\Sigma];\bC)$ by the full duals of $P_{(m+1)}\otimes 
\Omega^0[\Sigma]$ and $P_{(m+1)} \otimes \Omega^1[\Sigma]$, in cohomology 
degrees $0$ and $1$ and symmetric degrees $m+1$ and $m$, respectively. 
Generating co-cycles are constructed as in Theorem \ref{sym}, and the 
algebra is completed in the inverse limit topology defined by the 
order-of-pole filtration on $\Omega^i[\Sigma]$.
\end{theorem}

\section{Proof for truncated algebras}\label{2}
Assuming Theorem \ref{sym}, we now explain how Feigin's construction in 
\cite{feig2} proves Theorem \ref{trunc}, the conjecture for truncated Lie 
algebras. Its shadow is the specialisation $t=q^n$ in the combinatorial 
literature ($s=z^n$ in our notation). We can resolve $\frg[z]/z^n$ by the \textit{differential graded Lie algebra} $\left(\frg[z,s], 
\partial\right)$ with differential $\partial s = z^n$,
\begin{equation}\label{251}			
\left\{s\frg[z] \xrightarrow{\partial: s\mapsto z^n}\frg[z]\right\}
		\xrightarrow{\sim}\frg[z]/z^n.
\end{equation} 
This identifies $H^*\left(\frg[z]/z^n\right)$ with the \textit
{hyper-cohomology} of $\left(\frg[z,s],\partial\right)$, and $H^\bullet
\left(\frg[z]/z^n,\frg\right)$ with the relative one of the pair $\left(
(\frg[z,s],\partial),\frg\right)$. Recall that hyper-cohomology is computed 
by a \textit{double complex}, where Koszul's differential is supplemented 
by the one induced by $\partial$. This leads to a convergent spectral 
sequence, with
\begin{equation}\label{252}			
E_1^{p,q} = H_{res}^{q-p}\left(\frg[z],\frg; \sym^p(\frg[z])_{res}^*\right)
\Rightarrow H^{p+q}\left(\frg[z]/z^n,\frg\right). 
\end{equation}
The $E_1^{p,q}$ term arises by ignoring $\partial$, and is the portion of 
$H_{res}^{p+q}\left(\frg[z,s],\frg\right)$ with $s$-weight $(-p)$, cf.\ 
(\ref{18}). If we assign weight 1 to $z$ and weight $n$ to $s$, then
$\left(\frg[z,s],\partial\right)$ carries this additional \textit{$z$-grading},
preserved by $\partial$ and hence by the spectral sequence.

\begin{lemma}\label{253} 
Let $n>0$. $E_2^{p,q}$ is the free skew-commutative algebra generated by the 
dual of the sum of vector spaces 
$s^{\otimes m}\bC[z]dz/d\left(z^n\bC[z]\right)$,
placed in bi-degrees $(p,q)=(m,m+1)$, as $m$ ranges over $m_1,\ldots,m_\ell$.
The $z$-weight of $s$ is $n$.
\end{lemma}	
\begin{proof}[Proof of Theorem \ref{trunc}] The $E_2$ term of Lemma \ref{253} 
already meets the dimensional lower bound for our cohomology (Remark \ref
{12}.iii). Therefore, $E_2=E_\infty $ is the associated graded ring for a 
filtration on $H^\bullet\left(\frg[z]/z^n,\frg\right)$, compatible with the 
$z$-grading. However, freedom of $E_2$ forces $H^\bullet$ to be isomorphic 
to the same, and we get the desired description of $H^\bullet\left(\frg[z]/
z^n\right)$ from the factorisation (\ref{12}.i).
\end{proof}

\begin{proof}[Proof of Lemma \ref{253}] The description in Theorem \ref{sym} 
of the generating cocycles $E_\Phi$ and $S_\Phi$ of $E_1$ allow us to compute 
$\delta_1$. The $S_\Phi$ have nowhere to go, but for $E_\Phi:
\Lambda^1\otimes\sym^m\to\bC[z]dz$, we get
\begin{align}
\left(\delta_1 E_\Phi\right)(\sigma_0\cdot\ldots\cdot \sigma_m) 
&= E_\Phi\left(\partial(\sigma_0\cdot\ldots \cdot \sigma_m)\right)
	\nonumber\\
&=\sum\nolimits_k E_\Phi\left(z^n\sigma_k\otimes \sigma_0\cdot\ldots\cdot
\hat\sigma_k\cdot\ldots\cdot \sigma_m\right) 
	\nonumber\\
&=\sum\nolimits_k \Phi\left(\sigma_0\cdot\ldots\cdot d(z^n\sigma_k)
	\cdot\ldots\cdot\sigma_m\right) \\ 
&= (m+1)n\cdot z^{n-1}dz\cdot \Phi(\sigma_0\cdot\ldots\cdot\sigma_m) +
	z^n\cdot d\Phi\left(\sigma_0\cdot\ldots\cdot\sigma_m\right) 
	\nonumber\\ 
&=\left((m+1)n \cdot z^{n-1}dz\cdot S_\Phi + z^n\cdot dS_\Phi\right)
\left(\sigma_0\cdot\ldots\cdot\sigma_m\right), \nonumber
\end{align}
and so $\delta_1$ is the transpose of the linear operator $(m+1)n\cdot 
z^{n-1}dz\wedge +z^n\cdot d$, from $\bC[z]$ to $\bC[z]dz$. This has no kernel 
for $n>0$, and its co-kernel is $\bC[z]dz/d\left(z^n\bC[z]\right)$. 
\end{proof}

\begin{remark}\label{255}\begin{trivlist}\itemsep0ex  
\item(i) On $\frg[z,s]$, $\partial$ is given by the super-vector field 
$z^n\partial/\partial s$. This acts on the presentation (\ref{16}) of 
the homology primitives,
\begin{equation}\label{256}
z^n\partial/\partial s: \bC[z]\cdot s(ds)^m\to 
	\frac{\bC[z]\cdot (ds)^m + \bC[z]dz\cdot s(ds)^{m-1}}
	 {d\left(\bC[z]\cdot s(ds)^{m-1}\right)}. 
\end{equation}
Identifying the target space with $\bC[z]dz\cdot s(ds)^{m-1}$ by projection, 
we can check that $z^k\cdot s(ds)^m$ maps to $(mn+n+k)\cdot z^{n+k-1}dz
\cdot s(ds)^{m-1}$. This map agrees with (the dual of) the differential
$\delta_1$ in the preceding lemma, confirming our claim that the description 
(\ref{16}) was natural.
\item(ii) If $n=0$, the map in (\ref{256}) is surjective, with $1$-dimensional 
kernel; so $E_\infty^{p,q}$ now lives on the diagonal, and equals $(\sym^p
\frg^*)^\frg$. This is, in fact, a correct interpretation of $H^*(0,\frg;\bC)$. 
\end{trivlist}\end{remark}

\section{The Laplacian on the Koszul complex}\label{3}
In preparation for the proof of Theorem \ref{sym}, we now study the Koszul 
complex for the pair $(\frg[z,s],\frg)$ and establish the key formula \eqref
{311} for its Laplacian. 
 
\subsection{} For explicit work with $\frg[z]$-co-chains, we introduce the 
following derivations on $\Lambda\otimes\sym := \Lambda(\frg[z]/\frg)_{res}^*
\otimes\sym\frg[z]_{res}^*$, describing the brutally truncated adjoint action 
of $\frg[z,z^{-1}]$:
\begin{align}\label{32}		
\ad_a(m): \psi^b(n)\mapsto &\left\{ 
		\begin{array}{ll}\psi^{[a,b]}(m+n),&\mbox{if }m+n<0,\\
		0,&\mbox{if }m+n\ge 0;\end{array} \right.\\
R_a(m): \sigma^b(n)\mapsto &\left\{ 
		\begin{array}{ll}\sigma^{[a,b]}(m+n)&\mbox{if } m+n\le 0,\\ 
		0,&\mbox{if }m+n>0.\end{array} \right. \label{33}
\end{align}
Notations are as in \S\ref{02}, $m\in\bZ$ and $a,b$ range over $A:=\left\{1,
\ldots,\dim\frg\right\}$. Let 
\begin{equation}\label{34}
\bar\partial =\sum_{a\in A;m>0}\left\{\psi^a(-m)\otimes R_a(m) + 
	\psi^a(-m)\cdot \ad_a(m)\otimes 1/2\right\},
\end{equation} 
where $\psi^a(-m)$ doubles notationally for the appropriate multiplication
operator. The notation $\bar\partial$ stems from its geometric origin as
a Dolbeault operator on the loop Grassmannian of $G$.

\begin{definition}\label{35} 
The \emph{restricted Koszul complex} $\left(C^\bullet,\bar\partial\right)$ 
for the pair $(\frg[z],\frg)$ with coefficients in $\sym\frg[z]_{res}^*$ 
is the $\frg$-invariant part of $\Lambda^\bullet\otimes\sym$, with 
differential (\ref{34}).
\end{definition}

\subsection{The metric and the Laplacian.}\label{36} 
Define a hermitian metric on $\Lambda\otimes\sym$ by setting 
\[
\langle\sigma^a(m) | \sigma^b(n)\rangle = 1, 
\quad\langle \psi^a(m)|\psi^b(n)\rangle = -1/n, 
\quad\mbox{ if }m=n\mbox{ and }a=b,
\]
and both products to zero otherwise; we then take the multi-linear 
extension. Thus, $\sigma^a(m)^n/\sqrt{n!}$ has norm $1$. The hermitian 
adjoints to (\ref{32}) are the derivations defined by
\begin{equation}\label{32*}
\ad_a(m)^*\psi^b(n) = \frac{n-m}{n} \psi^{[a,b]}(n-m),\quad
\mbox{or zero, if } n\ge m.      \tag{\ref{32}$^*$}
\end{equation}
The $R$'s of (\ref{33}) satisfy the simpler relation $R_a(m)^*=R_a(-m)$. 
The adjoint of (\ref{34}) is 
\begin{equation}\label{34*}
\bar\partial^* = \sum_{a\in A;m>0} \left\{\psi^a(-m)^*\otimes R_a(-m) + 
\ad_a(m)^*\cdot \psi^a(-m)^*\otimes 1/2\right\}.     \tag{\ref{34}$^*$}
\end{equation}
A (restricted) Koszul cocycle in the kernel of the \textit{Laplacian} 
$\overline\square := \left(\bar\partial + \bar\partial^*\right)^2 = 
\bar\partial\bar\partial^* + \bar\partial^*\bar\partial$ is called \textit
{harmonic}. Since $\bar\partial$, $\bar\partial^*$ and $\overline\square$ 
preserve the orthogonal decomposition into the finite-dimensional $(z,s)
$-weight spaces, elementary linear algebra gives the following ``Hodge 
decomposition":

\begin{proposition}\label{37} 
The map from harmonic cocycles $\cH^k\subset C^k$ to their cohomology classes, 
via the decompositions $\ker\bar\partial = \im\bar\partial\oplus\cH^k$, 
$C^k = \im\bar\partial \oplus \cH^k \oplus \im\bar\partial^*$, is a linear 
isomorphism. \qed
\end{proposition}

To investigate $\overline\square$, we introduce the following adjoint pairs 
of operators: 
\begin{align}\label{38}
d_a(m)&:\sigma^b(n)\mapsto \psi^{[a,b]}(m+n),\mbox{ or zero, if }m+n\ge 0,
& d_a(m)\psi^b(n) &= 0;\\
d_a(m)^*&:\psi^b(n)\mapsto - \sigma^{[a,b]}(n-m)/n,\mbox{ or zero, if }n>m,
& d_a(m)^*\sigma^b(n) &= 0, \label{38*}\tag{\theequation $^*$}
\end{align}
extended to odd-degree derivations of $\Lambda\otimes\sym$. Finally, let	
\begin{align}\label{39}
D: &= \sum_{m>0;a\in A} d_a(-m)d_a(-m)^*, \\
\square: &= \sum_{a\in A;m>0}\frac{1}{m}\left[R_a(-m)+\ad_a(-m)\right]
	\left[R_a(m)+\ad_a(-m)^*\right].
\end{align}

\begin{theorem} \label{311}
On $C^\bullet$, $\overline\square=\square + D$. In particular, the harmonic 
forms are the joint kernel in $\Lambda \otimes \sym$ of the derivations 
$d_a(-m)^*$, as $a\in A$, $m>0$, and $R_a(m) + \ad_a(-m)^*$, as $a\in A$, 
$m\ge 0$.
\end{theorem}
\noindent It follows that the harmonic co-cycles form a sub-algebra, since 
they are cut out by derivations. We shall identify them in \S\ref{4}; the 
rest of this section is devoted to proving \eqref{311}. 

\begin{proof}[First proof of (\ref{311})] Introduce yet another operator 
\begin{equation}\label{312}
K:= \sum_{a,b\in A;\,m>0} \left(R_{[a,b]}(0) + \ad_{[a,b]}(0)\right)
	\cdot \psi^a(-m)\wedge\psi^b(-m)^*. 
\end{equation}
Note that the $\psi\wedge\psi^*$ factor could equally well be written 
in first position, because
\begin{multline*}
\sum_{a,b}\left[\ad_{[a,b]}(0),\psi^a(-m)\wedge\psi^b(-m)^*\right]\\
	=\sum_{a,b} \left(\psi^{[[a,b],a]}(-m)\wedge \psi^b(-m)^* 
		+ \psi^a(-m)\wedge \psi^{[[a,b],b]}(-m)^*\right) \\
=\sum_{a,b} \left(\psi^{[a,b]}(-m)\wedge\psi^{[a,b]}(-m)^*
		- \psi^{[a,b]}(-m)\wedge\psi^{[a,b]}(-m)^* \right) = 0.
\end{multline*}
As the first factor represents the total co-adjoint action of $\frg$ on 
$\Lambda\otimes\sym$, $K=0$ on the sub-complex $C^\bullet$ of $\frg
$-invariants, and Theorem \ref{311} is a special case of the following 
lemma.
\end{proof}

\begin{lemma}\label{313}
$\overline\square = \square + D + K$ on $\Lambda\otimes\sym$.
\end{lemma}
\begin{proof} All the terms are second-order differential operators on 
$\Lambda\otimes\sym$. It suffices, then, to verify the identity on 
quadratic germs. The brutal calculations are performed in the Appendix. 
\end{proof} 

\begin{proof}[Second proof of (\ref{311})] Let  $V$ be a negatively graded 
$\frg[z]$-module, such that $z^m\frg$ maps $V(n)$ to $V(n+m)$. Assume that 
$V$ carries a hermitian inner product, compatible with the hermitian 
involution on the zero-modes $\frg\subseteq \frg[z]$, and for which the 
graded pieces are mutually orthogonal. For us, $V$ will be $\sym\frg[z]_
{res}^*$. Write $R_a(m)$ for the action of $z^m\xi_a$ on $V$ and define, 
for $m\ge 0$, $R_a(-m):= R_a(m)^*$. Define $\square$ and $\overline\square$ 
as before; our conditions on $V$ ensure the finiteness of the sums. Define 
an endomorphism of $V\otimes\Lambda (\frg[z]/\frg)_{res}^*$ by the formula 

\begin{equation}\label{314} 
T_V^\Lambda:= \sum_{\ofrac{a,b\in A}{m,n>0}} 
		\left\{\left[R_a(m),R_b(-n)\right] - R_{[a,b]}(m-n) \right\}
				\otimes\psi^a(-m)\wedge \psi^b(-n)^*. 
\end{equation}
Our theorem now splits up into the two propositions that follow; the first 
is known as \textit{ Nakano's Identity}, the second describes $T_V^\Lambda$ 
when $V=\sym\frg[z]_{res}^*$. 
\end{proof}

\begin{proposition}[\cite{tel1}, Prop.~2.4.7]\label{315}
On $C^k$, $\overline\square = \square + T_S^\Lambda + k$.\qed
\end{proposition}

\begin{remark}\label{316}\begin{trivlist}\itemsep0ex  
\item(i) Our $R_a(m)$ is the $\theta_a(m)$ of \cite{tel1}, \S2.4, whereas
the operators $R_a(m)$ there are zero here, as is the level $h$. The constant
$2c$ from \cite{tel1} is replaced by $1$, because of our use of the Killing
form, instead of the basic inner product. A sign discrepancy in the definition
of $T_V^\Lambda $ arises, because our $\xi_a$ here are self-adjoint, and not
skew-adjoint as in \cite{tel1}.
\item(ii) \cite{tel1} assumed finite-dimensionality of $V$, but our grading 
condition is an adequate substitute. 
\end{trivlist}\end{remark}

\begin{proposition}\label{317}
On $\Lambda^k\otimes\sym$, $D = T_S^\Lambda + k$.
\end{proposition}
\subsection*{Proof.} Both sides are second-order differential operators on 
$\Lambda\otimes\sym$ and kill $1\otimes\sym$, so it suffices to check 
the equality on the following three terms of degree $\le 2$. Note that
$T_S^\Lambda = 0$ on $\Lambda\otimes 1$, and that $\sum\nolimits_a 
\psi^{[a,[a,b]]}(-n) = \psi^b(-n)$, because $\sum\nolimits_a \ad(\xi_a)^2 =
\mathbf{1}$ on $\frg$. 
\begin{align*}
D\psi^b(-n) &= \sum_{\ofrac{a\in A}{0<m\le n}} 
				d_a(-m)\sigma^{[a,b]}(m-n)/n\\ 
			&= \sum_{\ofrac{a\in A}{0<m\le n}}
				\psi^{[a,[a,b]]}(-n)/n = \psi^b(-n);
\end{align*}
\begin{align*}
D\left(\psi^b(-n)\wedge \psi^c(-p)\right) &= 
		D\psi^b(-n)\wedge \psi^c(-p) + \psi^b(-n)\wedge D\psi^c(-p)\\
	&= 2\cdot \psi^b(-n)\wedge\psi^c(-p);
\end{align*}
\begin{align*}
D\left(\sigma^c(-p)\cdot \psi^d(-q)\right) &= \sigma^c(-p)\cdot D\psi^d(-q)
	+ \frac{1}{q}\sum_{\ofrac{a\in A}{0<m\le q}}
		{\sigma^{[a,d]}(m-q)\cdot \psi^{[a,c]}(-m-p)}\\
		&= \sigma^c(-p)\cdot \psi^d(-q) + 
			T_S^\Lambda\left(\sigma^c(-p)\cdot \psi^d(-q)\right),
\end{align*}
with the last equality following from
\begin{multline*}
T_S^\Lambda \left(\sigma^c(-p)\cdot \psi^d(-q)\right)
	= \frac{1}{q} \sum_{\ofrac{a\in A}{0<m}}
		{\left\{\left[R_a(m), R_d(-q)\right] - R_{[a,d]}(m-q)\right\}
			\sigma^c(-p)\cdot \psi^a(-m)}\\
	= \frac{1}{q}\sum_{ \ofrac{a\in A}{0<m\le p+q}}
		{\sigma^{[a,[d,c]]}(m-q-p)\cdot \psi^a(-m)} - 
			\frac{1}{q}\sum_{\ofrac{a\in A}{0<m\le p}}
				\sigma ^{[d,[a,c]]}(m-q-p)\cdot \psi^a(-m) \\ 
	\qquad\qquad -\frac{1}{q}
		\sum_{ \ofrac{a\in A}{ {0<m\le p+q}}} 
			{\sigma^{[[a,d],c]}(m-q-p)\cdot \psi^a(-m)}\\ 
	= \frac{1}{q}\sum_{\ofrac{b\in A}{0<m\le p+q}}
		{\sigma^{[d,[a,c]]}(m-q-p)\cdot \psi^a(-m)} + \frac{1}{q}\sum_
			{ \ofrac{a\in A} {{0<m\le p}} }
			{\sigma^{[d,[a,c]]}(m-q-p)\cdot \psi^a(-m)}\\ 
	= \frac{1}{q}\sum_{ \ofrac{a\in A} { {p<m\le q}}}
			{\sigma^{[d,a]}(m-q-p)\cdot \psi^{[c,a]}(-m)}.\ \qed
\end{multline*}

\section{The harmonic forms and proof of Theorem \ref{sym}} \label{4}
We now use Theorem \ref{311} to identify the harmonic forms in $C^\bullet$; 
version (\ref{sym}) of the strong Macdonald conjecture follows by assembling 
Propositions \ref{45}, \ref{48} and \ref{410}.

\subsection{Relabelling $\psi$.} \label{41}
It will help to identify $\Lambda \left(\frg[z]/\frg\right)_{res}^*$ with 
$\Lambda\frg[z]_{res}^*$ by the isomorphism $d/dz: \frg[z]/\frg\cong \frg[z]$. 
This amounts to relabelling the exterior generators, with $\psi^a(-m)$ now 
denoting what used to be $(m+1)\cdot\psi^a(-m-1)$ ($m\ge 0$). Relations 
(\ref{32*}) and (\ref{38*}) now become
\begin{equation}\label{42}\begin{array}{rl}
\ad_a(-m)^*\psi^b(-n) =& \psi^{[a,b]}(m-n),\\
d_a(-m-1)^*\psi^b(-n) =& \sigma^{[a,b]}(m-n),\end{array}
\qquad\mbox{ or zero, if } m>n.
\end{equation}
According to (\ref{311}), the harmonic forms in the relative Koszul complex
(\ref{34}) are the forms in $\Lambda\frg[z]_{res}^*\otimes \sym\frg[z]
_{res}^*$ killed by $d_a(-m-1)^*$ and $R_a(m) + \ad_a(-m)^*$, as $m\ge 0$
and $a\in A$.

\subsection{The harmonic forms.} The graded vector space $\frg[[z,s]]:
=\frg[[z]]\oplus s\frg[[z]]$ carries the structure of a \textit{super-scheme}, 
if we declare functions to be the skew polynomial in finitely many of the 
components $z^m\frg$, $sz^m\frg$. It carries the adjoint action of the 
super-group scheme $G[[z,s]]$, which is a semi-direct product $G[[z,s]]
\cong G[[z]]\ltimes s\frg[[z]]$. 

\begin{lemma}\label{44} 
Identifying $\Lambda\frg[z]_{res}^*\otimes\sym\frg[z]_{res}^*$ with the 
(skew) polynomials on $\frg[[z,s]]$, the operators $d_a(-m-1)^*$ and 
$R_a(m) + \ad_a(-m)^*$, as $m\ge 0$, generate the co-adjoint action 
of $\frg[z,s]$.
\end{lemma}
\begin{proof} This is clear from (\ref{42}): $d_a(-m-1)^*$ is 
the co-adjoint action of $s\cdot z^m\xi_a$.
\end{proof} 

\begin{proposition}\label{45} The harmonic forms in $C^\bullet $ 
correspond to those skew polynomials on $\frg[[z,s]]$ which are invariant 
under the adjoint action of $G[[z,s]]$.
\end{proposition}
\begin{proof} Lie algebra and group invariance of functions are equivalent, 
because the action is locally-finite and factors, locally, through the 
finite-dimensional quotients $\frg[z,s]/z^N$.
\end{proof}

\begin{remark}\label{46}
The super-language can be avoided by identifying $\frg[[z,s]]$ with the 
tangent bundle to its even part $\frg[[z]]$, having declared the tangent 
spaces to be odd: the skew polynomials become the \textit{polynomial 
differential forms} on $\frg[[z]]$, and the invariant skew functions under 
$G[[z,s]]$ correspond to the \textit{basic} forms under the Ad-action 
of $G[[z]]$. 
\end{remark}

\subsection{The invariant skew polynomials.} The (GIT) quotient $\frg/\!/G 
:= \mathrm{Spec}(\sym\frg^*)^G$ is the space $P$ of primitives in the 
co-algebra $\sym\frg/[\frg,\sym\frg]$. The quotient map $q:\frg\to P$ 
induces a morphism $Q:\frg[[z,s]]\to P[[z,s]]$, which is invariant under 
the adjoint action of $G[[z,s]]$.
\begin{proposition}\label{48}
The $\ad$-invariant skew polynomials on $\frg[[z,s]]$ are precisely the 
pull-backs by $Q$ of the skew polynomials on $P[[z,s]]$.
\end{proposition}
\begin{proof} Elements $\Lambda\frg[z]_{res}^*\otimes\sym\frg[z]_{res}^*$ 
are algebraic sections of the vector bundle $\Lambda\frg[z]_{res}^*$ over 
$\frg[[z]]$. As such, they are uniquely determined by their restriction to 
Zariski open subsets. The analogue holds for $P$. Now, the open subset $\frg
^{rs}\subset \frg$ of \textit{regular semi-simple} elements is an algebraic 
fibre bundle, via $q$, over the open subset $P^r\subset P$ of regular 
conjugacy classes. Let $\frg^{rs}[[z,s]]$ be the pull-back of $\frg^{rs}$ 
under the evaluation morphism $s=z=0$. Because of the local product 
structure, it is clear that $\ad$-invariant polynomials over $\frg^{rs}
[[z,s]]$ are precisely the pull-backs by $Q$ of functions on $P^r[[z,s]]$. 
In particular, the pull-back of polynomials from $P[[z,s]]$ to $\frg[[z,s]]$ 
is injective.

Let now $f$ be an invariant polynomial on $\frg[[z,s]]$. Its restriction 
to $\frg^{rs}[[z,s]]$ has the form $g\circ Q$, for some regular function 
$g$ on $P^r[[z,s]]$. Let $\frg^r\subset \frg$ be the open subset of 
\textit{regular} elements. A theorem of Kostant's ensures that $q: \frg^r
\to P$ is a submersion. In particular, it has local sections everywhere, 
so the morphism $Q:\frg^r[[z,s]]\to P[[z,s]]$ has local sections also. 
We can use local sections to extend our $g$ from $P^r[[z,s]]$ to $P[[z,s]]$,
because $f$ was everywhere defined upstairs. The extension of $g$ is unique,
and its $Q$-lifting must agree with $f$ everywhere, as it does so on an
open set. So we have written $g$ as a pull-back.
\end{proof}
 
\subsection{Relation to $S_\Phi$ and $E_\Phi$.} A polynomial $\Phi$ on 
$P$ defines a map $P[[z]]\to \bC[[z]]$ by point-wise evaluation, and the 
$m$th coefficient $\Phi(-m)$ of the image series is a polynomial on 
$P[[z]]$. The analogue holds for differential forms, or skew polynomials 
on our super-schemes \eqref{46}.

\begin{proposition}\label{410}
Let $\Phi_1,\ldots,\Phi_\ell$ be a basis of linear functions on P and let 
$\Phi_k(m)$ ($m\le 0$) be the associated Fourier mode basis of linear 
functions on $P[[z]]$. After $\psi$-relabelling as in \S\ref{41}, the 
cocycles $S_k(m)$ and $E_k(m)$ associated to $\Phi_k$ in (\ref{sym}) are 
the $Q$-lifts of $\Phi_k(m)$ and $d\Phi_k(m)$. 
\end{proposition}
\begin{proof} For $S_k(m)$, this is the obvious equality $\Phi_k(m)\circ Q = 
(\Phi_k\circ q)(m)$, the $(-m)$th Fourier mode of $\Phi\circ q$ on $\frg[[z]]$. 
For $E_k(m)$, observe that when replacing skew polynomials on $X[[z,s]]$ 
by forms on $X[[z]]$ as in Remark \ref{46} ($X=\frg,P$), $Q$ is the 
differential of its restriction $\frg[[z]]\to P[[z]]$, while $E_k(m-1)= 
dS_k(m)$, after our relabelling. 
\end{proof}

\subsection{The super-Iwahori algebra.} \label{411}
We now deduce Theorem \ref{111} from \ref{sym}. Let $\exp(\frB)$ be the 
closed Iwahori subgroup of $G[[z]]$, whose Lie algebra is the $z$-adic 
completion $\frB_z$ of $\frB$. We write $H_{\exp(\frB)}^\bullet(V)$, 
$H_{G[[z]]}^\bullet(V)$ for the \textit{algebraic group cohomologies} of 
$\exp(\frB)$, resp.\ $G[[z]]$ with coefficients in a representation $V$. 
Applying van Est's spectral sequence gives 
\begin{align*}		
H^\bullet\left(\frB,\frh;\sym\frB_{res}^*\right) 
		&= H_{\exp(\frB)}^\bullet\left(\sym\frB_{res}^*\right),\\	
H_{res}^\bullet\left(\frg[z],\frg;\sym\frg[[z]]_{res}^*\right) 
		&= H_{G[[z]]}^\bullet\left(\sym\frg[z]_{res}^*\right). 
\end{align*}				
We now relate the right-hand terms using \textit{Shapiro's spectral 
sequence}
\[
E_2^{p,q} = H_{G[[z]]}^p\left(R^q\,\mathrm{Ind}_{\exp(\frB)}^{G[[z]]} 
		\sym\frB_{res}^*\right)\Rightarrow H_{\exp(\frB)}^{p+q}
				\left(\sym\frB_{res}^*\right),
\]
whose collapse is a consequence of the following lemma, which, combined 
with the freedom of $\sym\frh^*$ as a $(\sym\frg^*)^\frg$-module, also 
completes the proof of Theorem \ref{111}. Write $R^q\,\mathrm{Ind}$ 
for $R^q\,\mathrm{Ind}_{\exp(\frB)}^{G[[z]]}$.

\begin{lemma}\label{412} 
$\mathrm{Ind}\;\sym\frB_{res}^* = \sym\frg[z]_{res}^*\otimes_
{(\sym\frg^*)^\frg}\sym\frh^*$, with the adjoint action of $G[[z]]$ 
on the first factor on the right; whereas $R^q\,\mathrm{Ind}\; 
\sym\frB_{res}^* =0 $ for $q>0$. 
\end{lemma}

\begin{proof} $R^q\,\mathrm{Ind}\;\sym\frB_{res}^*$ is the $q$th sheaf 
cohomology of the algebraic vector bundle $\sym\frB_{res}^*$ over the 
quotient variety $G[[z]]/\exp(\frB)\cong G/B$, and hence also the 
$q$th cohomology of the structure sheaf $\cO$ over the variety $G[[z]] 
\times_{\exp(\frB)}\frB_z$, with the adjoint action of $\exp(\frB)$ on 
$\frB_{z}$. Splitting $\frB_z$ as $\frb\times z\frg[[z]]$ and shearing 
off the second factor identifies this variety with $(G\times_B\frb)\times 
z\frg[[z]]$. The factor $G\times_B\frb$ maps properly and generically 
finitely to $\frg$ via $\mu:(g,\beta)\mapsto g\beta g^{-1}$. The canonical 
bundle upstairs is trivial, and a theorem of Grauert and Riemenschneider 
ensures the vanishing of higher cohomology of $\cO$, and thus of the higher 
$R^q\,\mathrm{Ind}$'s. 

The functions on $G\times_B\frb$ are identified with $\sym\frh^*\otimes 
_{(\sym\frg^*)^G}\sym\frg^*$ by the Stein factorisation of $\mu$, 
\[
G\times_B\frb\xrightarrow{(\pi,\mu)}\frh\times_{\frg/\!/G}\frg\to \frg, 
\]
where $\pi:\frb\to\frh$ is the natural projection and the second arrow the 
second projection. (The middle space is regular in co-dimension three, 
therefore normal.) Using this and evaluation at $z=0$, we can factor the 
conjugation morphism $G[[z]]\times_{\exp(\frB)}\frB_z \to\frg[[z]]$ into 
the $G[[z]]$-equivariant maps below, of which the first has proper and 
connected fibres,
\[
G[[z]]\times_{\exp(\frB)}\frB_z\to 
		\frh\times_{\frg/\!/G}\frg[[z]]\to \frg[[z]].
\]
This exhibits the space of functions $\mathrm{Ind}\;\sym\frB_{res}^*$ 
on $G[[z]]\times_{\exp(\frB)}\frB_z$ to be as claimed. 
\end{proof}

\part{Hodge theory}
We now turn to a remarkable application of the strong Macdonald theorem: 
the determination of Dolbeault cohomologies $H^q(\Omega^p)$ and the 
Hodge-de Rham sequence for flag varieties of loop groups. For the loop 
Grassmannian $X$, these are described formally from $H^\bullet(BG)$ 
and de Rham's operator $d:\bC[[z]]\to\bC[[z]]dz$ on the formal disk 
(Theorem \ref{thinhodge}). In particular, we find that the sequence 
collapses at $E_2$, and not at $E_1$, as in the case of smooth projective 
varieties. This failure of Hodge decomposition is unexpected, given the 
(ind-)projective nature of $X$; surprisingly for a homogeneous space, 
the explanation lies in the lack of smoothness. 
  
Similar results hold for other flag varieties, associated as in \S\ref{7} 
to a smooth affine curve $\Sigma$; the Dolbeault groups and first 
differentials in the Hodge sequence arise from $d: \Omega^0[\Sigma] 
\to\Omega^1[\Sigma]$ (Theorem \ref{thickhodge}). This concords 
with the Hodge decomposition established in \cite{tel4} for the \textit
{closed curve} analogue of our flag varieties, the moduli stack of $G
$-bundles over a smooth projective curve. Evidently, the failure of Hodge 
decomposition for flag varieties is rooted in the same phenomenon 
for open curves, but we do not feel that we have a satisfactory 
explanation.  

The description of Dolbeault groups is unified conceptually in \S\ref
{unified}, where we construct generating co-cycles. We also interpret 
the Macdonald cohomology of Chapter I as the Dolbeault cohomology of 
the \textit{classifying stack} $BG[[z]]$. That is also the moduli stack 
of principal $G$-bundles on the formal disk; its relevance arises by 
viewing the flag varieties as moduli space of $G$-bundles over the 
completion of $\Sigma$, trivialised in a formal neighbourhood of the 
divisor at infinity. The construction leads to the proofs in \S\ref
{proof}, and our arguments feed back in \S\ref{10} into some new Lie 
algebra results, including the proof of Theorem \ref{1-15} on the 
cohomology of $\frg[\Sigma,s]$.

To keep the statements straightforward, $G$ will be \textit{simple} and 
\textit{simply connected}. 

\section{Dolbeault cohomology of the loop Grassmannian}
\label{5}
\subsection{The loop Grassmannian.}\label{61} 
By the \textit{loop group $LG$} of $G$ we mean the group $G((z))$ of formal 
Laurent loops; it is an \textit{ind-group-scheme}, filtered by the order 
of the pole. (The order, but not the ind-structure, depend on a choice of 
closed embedding $G$ into affine space.) The \textit{loop Grassmannian} of 
$G$ is the quotient (ind-)variety $X:=LG/G[[z]]$ of $LG$. This is \textit
{ind-projective}---an increasing union of closed projective varieties---and 
in fact Kodaira-embeds in a direct limit projective space \cite{kum}. The 
largest ind-projective quotient of $LG$ is the \textit{full flag variety} 
$LG/\exp(\frB)$, which is a bundle over $X$ with fibre the full $G$-flag 
manifold $G/B$; the other ind-projective quotients correspond to the 
subgroups of $LG$ containing $\exp(\frB)$.    

As a homogeneous space, $X$ is formally smooth, so there is an obvious 
meaning for the algebraic differentials $\Omega^p$. The Dolbeault 
cohomologies\footnote{We retain the analytic term Dolbeault cohomology to 
indicate the presence of differential forms, even when using \textit
{algebraic} sheaf cohomology; the distinction is immaterial for $X$, by 
GAGA.} $H^q(X;\Omega^p)$ carry a translation action of the loop group, and 
a grading from the $\bC^\times$-action scaling $z$ (loop rotation). 
\begin{proposition}\label{62}
$H^\bullet(X; \Omega^\bullet)$ is the direct product of its $z$-weight 
spaces, and the action of $LG$ is trivial.
\end{proposition}  
\begin{proof}
The sheaves $\Omega^p$ are associated to the co-adjoint action of $G[[z]]$ 
on the full duals of the exterior powers of $\frg((z))/\frg[[z]]$; as such, 
they carry a decreasing filtration $Z^n\Omega^p_{n>0}$ by $z$-weight, and 
are complete thereunder. The associated sheaves $\Gr^n\Omega^p$ are sections 
of finite-dimensional bundles, stemming from the co-adjoint action of 
$G[[z]]$ on $\Gr^n\Lambda^p \left\{\frg((z))/ \frg[[z]]\right\}^*$. This 
factors through $G$ by the evaluation $z=0$. The cohomologies of the 
$\Gr^n\Omega^p$ are then finite-dimensional, trivial $LG$-representations 
\cite{kum}; so, then, are the cohomologies $H^*(X;\Omega^p/Z^n\Omega^p)$ 
of the $z$-truncations, which are finite extensions of such representations. 

The $\Omega^p/Z^n\Omega^p$ give a surjective system of sections over any 
ind-affine open subset of $X$. The Mittag-Leffler condition for their 
cohomologies is clear by finite-dimensionality; we conclude the equality 
$H^*(X;\Omega^p)=\lim_n H^*(X;\Omega^p/Z^n\Omega^p)$ and the proposition. 
\end{proof}

Our main theorem describes the Dolbeault groups of $X$ and the action 
thereon of de Rham's operator $\partial:\Omega^p\to\Omega^{p+1}$. The  
$z$-adic completeness, ensured by the previous proposition, stems from 
the close relation of $X$ with the formal disk (cf.\ the discussion of 
thick flag varieties in \S\ref{7}). 

\begin{maintheorem}\label{thinhodge}\hfill\begin{enumerate}
\item $H^\bullet(X;\Omega^\bullet)$ is the $z$-adically completed skew 
power series ring generated by copies of $\bC[[z]]$ and $\bC[[z]]dz$, 
lying in $H^m(\Omega^m)$ and $H^m(\Omega^{m+1})$, respectively ($m = 
m_1,\ldots,m_\ell$). 
\item De Rham's differential $\partial:H^q(X;\Omega^p)\to H^q(X; \Omega^
{p+1})$ is the derivation induced by $d:\bC[[z]]\to\bC[[z]]dz$ on generators. 
Its cohomology is the free algebra on $\ell$ generators in bi-degrees $(m,m)$, 
ranging over the exponents of $\frg$.
\end{enumerate}
\end{maintheorem}
\noindent The generators are constructed in Proposition \ref{86}, and the 
theorem will be proved in \S\ref{proof}.  

\subsection{Failure of Hodge decomposition.} \label{failhodge}
In the analytic topology, de Rham's complex $(\Omega^\bullet,\partial)$ 
resolves the constant sheaf $\bC$. GAGA implies that the hyper-cohomology 
$\bH^\bullet(X;\Omega^\bullet,\partial)$ agrees with the complex cohomology 
$H^\bullet(X;\bC)$. Recall \cite{gar} that $X$ is homotopy equivalent to the 
group $\Omega G$ of based continuous loops, or again, to the double loop 
space $\Omega^2 BG$ of the classifying space. Its complex cohomology is 
freely generated by the $S^2$-transgressions of the generators of $H^\bullet
(BG;\bC)\cong(\sym^{2\bullet}\frg^*)^G$. Theorem \ref{thinhodge} implies 
that the differential $\partial_1$ on $H^q(\Omega^p)$ resolves the complex 
cohomology of $X$. In other words, the \textit{Hodge-de Rham spectral 
sequence} induced by $\partial$ on $\Omega^\bullet$ collapses at $E_2$. 

As $X$ is ind-projective, formally smooth and reduced \cite{ls}, we might 
have expected a Hodge decomposition of its complex cohomology into the 
$H^q(X;\Omega^p)$. Failure of this is has the following consequence, 
as announced in \cite{tel2}. The proof is lifted from \cite{simtel}, \S7. 
We emphasise that the result asserts more than the absence of a global
expression for $X$ as a union of smooth projective subvarieties (indeed, 
there is a cleaner argument for this last fact, \cite{gro}).  

\begin{theorem}\label{65}
$X$ is not a smooth complex manifold: that is, it cannot be expressed, 
locally in the analytic topology, as an increasing union of smooth 
sub-manifolds.
\end{theorem}

\noindent Because $X$ is homogeneous, it is singular everywhere. The same
will be true for the full flag variety $LG/\exp(\frB)$, and for the loop 
group $LG$ itself.
\begin{proof}[Proof of (\ref{65})]
Expressing $X$ as a union of projective sub-varieties $Y_n$ (for instance,
the closed Bruhat varieties) gives an equivalence of $X$ with the ($0$-stack) 
represented, over the category of complex schemes of finite type, by the 
groupoid $\coprod Y_n\rightrightarrows \coprod Y_n$. (The two structural 
maps are the identity and the family of inclusions $Y_n\hookrightarrow Y_
{n+1}$.) In more traditional terms, this gives a \textit{simplicial 
resolution} $Y_\bullet\xrightarrow {\vep} X$ by a simplicial variety whose 
space of $n$-simplices is a union of projective varieties, for each $n$. 
Resolution of singularities and the method of hyper-coverings in \cite{del} 
allows us to replace $Y_\bullet$ by a smooth simplicial resolution $X_\bullet
\xrightarrow{\vep} X$ (in the topology of generated by proper surjective 
covers). The total direct image $R\vep_*$ of de Rham's complex 
$(\Omega^\bullet,\partial; F)$ with the Hodge filtration 
\[
F^p\Omega^\bullet:= \left[\Omega^p\xrightarrow{\partial}\Omega^{p+1}
		\xrightarrow{\partial}\ldots\right]
\]
is the \textit{DuBois complex}  \cite{dub} on $X$. The associated graded 
complex $\underline{\Omega}^p:= \Gr^p R\vep_*(\Omega^\bullet,\partial; F)$
is the `correct' singular-variety analogue of the $p$th Hodge-graded 
sheaf of the constant sheaf $\bC$. Because $X_\bullet$ is simplicially 
projective, the cohomology of $\underline{\Omega}^p$ satisfies the 
Hodge decomoposition
\begin{equation}\label{66}
H^n(X;\bC)\cong \bigoplus_{p+q=n}H^q(X;\underline{\Omega}^p).
\end{equation} 
The key properties of the DuBois complex are \textit{locality in the 
analytic topology} and \textit{independence of simplicial resolution}.
The restriction in \cite{dub} to finite-dimensional varieties need not 
trouble us: the arguments there show that $\underline{\Omega}^p$ is a 
well-defined, up to canonical isomorphism, in the bounded-below derived 
category of coherent sheaves over the site of analytic spaces, in the 
topology generated by both proper surjective maps and open covers. We 
are studying the hyper-cohomology of these $\underline{\Omega}^p$ in the 
restricted site of spaces over $X$. These properties would lead to a 
quasi-isomorphism $\underline{\Omega}^p \sim \Omega^p$, if $X$ was a
complex manifold in the sense of Theorem \ref{65}. But then, \eqref{66} conflicts with Theorem \ref{thinhodge}.  
\end{proof}


\section{Application: a $_1\psi_1$ summation}\label{6} 
The $H^q(X;\Omega^p)$ are graded by $z$-weight, and the $z$-weighted 
holomorphic Euler characteristics, for various $p$, can be collected in 
the \textit{$E$-series} 
\begin{equation}\label{68}
E(z,t):= \sum\nolimits_{p,q} (-1)^q (-t)^p 
		\dim_z H^q\left(X,\Omega^p\right)\in \bZ[[z,t]].
\end{equation}
\vskip-1ex
\subsection{Kac formula.}\label{euler}
The Mittag-Leffler conditions in the proof of Proposition \ref{62} 
imply the convergence of the spectral sequence for the $Z$-filtration,
\[
E_1^{r,s} = H^{r+s}\left(X; \Gr^r\Omega^p\right) \Rightarrow 
		H^{r+s}\left(X; \Omega^p\right),
\] 
whence it follows that the characteristic is already computed by $E_1$. 
Because $\Gr\,\Omega^p$ is a product of bundles associated to irreducible 
representations of $G[[z]]$, the $z$-weighted Euler characteristics are 
described explicitly by the \textit{Kac character formula} \cite{kac}. 
Choose a maximal torus $T\subset G$ and recall that the \textit{affine 
Weyl group} $\waff$ is the semi-direct product of the finite Weyl group 
by the co-root lattice; it acts on Fourier polynomials on $T$ and in $z$, 
whereby a co-root $\gamma$ sends the Fourier mode $\me^\lambda$ of $T$ 
to $z^{\langle\lambda|\gamma\rangle} \me^\lambda$. (The Weyl group acts 
in the obvious way, and $z$ is unaffected.) The desired formula is 
the infinite sum of infinite products indexed by affine roots,
\begin{equation}\label{kac}
\sum_{w\in\waff} \prod_{\textstyle\ofrac{n>0}{\alpha}}
		w\left(\frac{1-tz^n\me^\alpha}{1-z^n\me^\alpha}\right)
		\cdot\prod_{\alpha>0}w(1-\me^\alpha)^{-1}
		\cdot\prod_{n>0}\left[\frac{1-tz^n}{1-z^n}\right]^\ell. 
\end{equation}
The summands are the $w$-transforms of the quotient of the $(T,z,t)
$-character of the fibre $\sum_p(-t)^p\Gr\,\Omega^p$ at the base-point 
of $X$ by the Kac denominator. The sum expands into a formal power series 
in $z$ and $t$, with characters of $T$ as coefficients.

\subsection{Relation to Ramanujan's $_1\psi_1$ sum.}\label{1psi1}
Factoring affine Weyl elements as $\gamma\cdot w$ (co-root times finite 
Weyl element) and leaving out, for now, the third factor converts 
\eqref{kac} into 
\[
\sum_{\gamma} \prod_{\textstyle\ofrac{n>0}{\alpha}}
		\frac{1-tz^{n+\langle \alpha | \gamma\rangle}\me^\alpha}
			{1-z^{n+\langle \alpha | \gamma\rangle}\me^\alpha}
		\cdot\sum_{w\in W}\prod_{\alpha>0}
			(1-z^{\langle w\alpha |\gamma\rangle}\me^{w\alpha})^{-1},
\]
having substituted $\alpha\mapsto w\alpha$ in the first product, to make 
it $w$-independent. The second factor, the sum over $W$, is identically $1$, 
by the Weyl denominator formula. Therefore, equating \eqref{kac} with our 
answer in Theorem \ref{thinhodge} gives the following identity:
\[
\sum_{\gamma} \prod_{\textstyle\ofrac{n>0}{\alpha}}
		\frac{1-tz^{n+\langle\alpha | \gamma\rangle} \me^\alpha}
			{1-z^{n+\langle\alpha | \gamma\rangle}\me^\alpha} =
		\prod_{\textstyle\ofrac{1\le k\le \ell}{n\ge 0}}
				\frac{(1-z^{n+1})(1-t^{m_k+1}z^{n+1})}
				{(1-tz^{n+1})(1-t^{m_k}z^n)}.
\]
It is part of the statement that the left-hand side is constant, as a 
function on $T$. 

For $G=\SL_2$, we obtain after setting $\me^\alpha =u$ the identity
\[
\sum_{m} \prod_{n>0}\frac{(1-tz^{n+2m} u^2)(1-tz^{n-2m} u^{-2})}
			{(1-z^{n+2m} u^2)(1-z^{n-2m} u^{-2})} =
		\frac{1}{1-t}\prod_{n>0}
				\frac{(1-z^n)(1-t^2 z^n)}
				{(1-tz^n)^2}
\]
which also follows from a 3-variable specialisation of Ramanujan's $_1\psi
_1$ sum (\cite{tel2}, \S5). (Note that our sum contains the even terms 
only; the ``other half" of the specialised $_1\psi_1$ sum is carried by 
the twisted $\SL_2$ loop Grassmannian, the odd component of $LG/G[[z]]$ 
for $G=\mathrm{PSL}_2$.) Thus, Theorem \ref{thinhodge} is a \textit{strong 
form} of (specialised) $_1\psi_1$ summation, generalised to (untwisted) 
affine root systems. We later learnt that (the ``weak" forms of) such 
generalised summation formulae, for all affine root systems, were 
independently discovered and proved by Macdonald \cite{mac}.

\section{Thick flag varieties}\label{7}
Related and, in a sense, opposite to $X$ is the quotient variety $\bfX
:= LG/G[z^{-1}]$. This is a scheme covered by translates of the open 
cell $\bfU \cong G[[z]]/G$, the $G[[z]]$-orbit of $1$. Generalisations of 
$\bfX$ are associated to smooth affine curves $\Sigma$, with divisor at 
infinity $D$ in their smooth completion $\oSig$. They are the quotients 
$\bfX_\Sigma:= L^DG/G[\Sigma]$ of a product $L^DG$ of loop groups, 
defined by local coordinates centred at the points of $D$, by the 
ind-subgroup $G[\Sigma]$ of $G$-valued regular maps. Variations decorated 
by bundles of $G$-flag varieties, attached to points of $\Sigma$, also 
exist, and our results extend easily to those, although we shall not 
spell that out. When a distinction is needed, we call the $\bfX_\Sigma$ 
and their variations \textit{thick flag varieties} of $LG$. 

\subsection{Relation to moduli spaces.} \label{71} One formulation of the 
\textit{uniformisation theorem} of \cite{ls} equates $X_\Sigma$ with the 
moduli space pairs $(\cP,\sigma)$ of algebraic principal $G$-bundles $\cP$ 
over $\oSig$, equipped with a section $\sigma$ over the formal neighbourhood 
$\wD$ of the divisor at infinity. In other words, $\bfX_\Sigma$ is the 
moduli space of relative $G$-bundles over the pair $(\oSig, \wD)$, and we 
also denote it by $\frM(\oSig,\wD)$. $\frM$ stands here for the \textit
{stack of morphisms to $BG$}, the classifying stack of $G$ (Appendix B of 
\cite{tel3}); thus, $\frM(\overline\Sigma)$ is the moduli stack of $G$-bundles over the closed curve. The corresponding description of $X$ is the moduli 
space of pairs, consisting of a $G$-bundle over $\bP^1$ and a section over 
$\bP^1 \setminus\{0\}$; this is the moduli space $\frM(\bP^1,\bP^1 \setminus
\{0\})$ of bundles over the respective pair. In this sense, $X$ is the $\bfX$ 
associated to the formal disk around $0$. Slightly more generally, $\frM 
(\oSig,\Sigma)$ is the product of loop Grassmannians associated to the 
points of $D$.

The thick flag varieties are smooth in an obvious geometric sense: the 
open cell in $\bfX$ is isomorphic to the vector space $\frg[[z]]/\frg$, 
while $\bfX$ is a principal $G[[z]]$-bundle over $\frM(\oSig)$. In 
their case, failure of Hodge decomposition should be attributed to 
``non-compactness".

\subsection{Technical note on spaces.} \label{technote}
We shall use the terms \textit{space} or, abusively, \textit{variety}, for 
the homogeneous spaces of $LG$. They live in a suitable world of contravariant 
functors on complex schemes: thus, the functor $\bfX_\Sigma$ sends a scheme 
$S$ to the set ``$\mathrm{Hom}\,(S,\bfX_\Sigma)$" of isomorphism classes of 
bundles over $(S\times\oSig, S\times\wD)$, and the ambient world is the 
category of sheaves over the topos of complex schemes of finite type, in the 
smooth (or \'etale) topology.\footnote{To include stacks, we must enrich the 
structure to include the simplicial sheaves and their homotopy category; see 
\cite{tel3} for a working introduction to this jargon.} 

The expert may raise some concern about the dual nature of thick flag
varieties: as they are schemes of \textit{infinite type} in their own 
right, the correct functorial perspective places them in the topos of 
\textit{all} complex schemes. Restricting to schemes of finite type 
converts these infinite-type schemes to the associated \textit{pro-objects}, 
with respect to the obvious filtration of $G[[z]]$ induced by the powers of 
$z$. However, these fibre, in (infinite) affine spaces, over smooth schemes 
of finite type (moduli of bundles over $\oSig$ with level structure). The 
results we use below (de Rham's theorem, vanishing of higher coherent 
cohomology) hold both for infinite affine space $\mathrm{Spec}\,\bC[z_1,
z_2,\ldots]$ and for the associated pro-affine space. Thus, our restriction 
to finite type is innocuous.\footnote{Restriction to finite type is only 
truly used for the thin varieties, in relation to the Du Bois complex; the 
reader is free to treat the thick flag varieties as schemes throughout.}

\subsection{Cohomology and Hodge structure.} \label{cohodge}
Recall now the analogue of the homotopy equivalence $X\sim\Omega G$ for 
thick varieties $\bfX_\Sigma$. The natural morphism from $\bfX_\Sigma = 
\frM(\oSig, \wD)$ to the stack $\frM(\oSig,D)$ of $G$-bundles on $(\oSig, D)$ 
(trivialised over $D$) is a fibre bundle in affine spaces; in particular, it 
is a homotopy equivalence. Similarly to Theorem 1' of \cite{tel3} (in which 
$D=\emptyset$), this last stack has the homotopy type of the space of the 
continuous maps from $\oSig$ to $BG$, based at $D$; the equivalence is the 
forgetful functor from the stack of ($D$-based) analytic bundles to that of 
continuous bundles.\footnote{This can be seen from the Atiyah-Bott 
construction of $\frM(\oSig)$.}

Generators of the algebra $H^\bullet(\frM(\oSig,D),\bQ)$ arise by 
transgressing those of $H^\bullet(BG)$ along a basis of cycles in 
$H_\bullet(\oSig,D)$; the latter is also the \textit{Borel-Moore homology} 
$H_\bullet^{BM} (\Sigma)$. As the classifying morphism $(\oSig,D)\times
\frM(\oSig,D)\to BG$ for the universal bundle is algebraic, the construction 
of generating classes is compatible with Hodge structures and we obtain 
(cf.~\cite{tel3}, Ch.~IV)

\begin{proposition}\label{cohom}
$H^\bullet(\frM(\oSig,D))$, with its Hodge structure, is the free algebra 
generated by $\prim\, H^\bullet(BG)\otimes H_\bullet^{BM}(\Sigma)$, with 
the natural Hodge structures on the factors.\qed 
\end{proposition} 

\noindent Recall \cite{del} that the Hodge structure on $BG$ is pure of 
type $(p,p)$. We can use the isomorphism $H^\bullet(\bfX_\Sigma)\cong 
H^\bullet(\frM(\oSig,D))$ to define the Hodge structure on $\bfX_\Sigma$, 
which is a scheme of infinite type. (By the argument in \S\ref{technote}, 
it agrees with the structure of the functor represented by $\bfX_\Sigma$ 
over the schemes of finite type). 

\subsection{Differentials.} Denote by $\Omega^p$ the sheaf of algebraic 
differential $p$-forms on any of our flag varieties. On $X$, this is the 
sheaf of sections of a pro-vector bundle, dual to $\Lambda^p TX$, but on 
thick flag varieties, it corresponds to an honest vector bundle, albeit 
of infinite rank. The is a de Rham differential $\partial: \Omega^p 
\to\Omega^{p+1}$. 
\begin{proposition}[Algebraic de Rham] 
$\bH^\bullet\left(\bfX_\Sigma;(\Omega^\bullet,\partial)\right)= H^\bullet
(\bfX_\Sigma;\bC)$, the former being the algebraic sheaf (hyper)cohomology, 
the latter defined in the analytic topology. 
\end{proposition}
\begin{proof} 
For $\bfX$, we use the standard \v{C}ech argument for the covering by the 
affine Weyl translates of the open cell; each finite intersection of the 
covering sets is a complement of finitely many coordinate hyperplanes in 
$\frg[[z]]/\frg$, where de Rham's theorem is obvious. The more general 
$\bfX_\Sigma$ are bundles in affine spaces over the (smooth, locally Artin) 
stacks $\frM(\oSig,D)$; de Rham's theorem for the total space follows from 
its knowledge on the fibres and on the base. 
\end{proof}

There results a convergent Hodge-de Rham spectral sequence 
\begin{equation}\label{hdR}
E_1^{p,q} = H^q(\bfX;\Omega^p),\qquad E_\infty^{p,q} 
		= \mathrm{Gr}^p H^{p+q}(\bfX;\bC),
\end{equation}
with the graded parts $\mathrm{Gr}^p$ of $H^*$ associated to the \textit
{na\"ive Hodge filtration}, the images of the truncated hyper-cohomologies 
$\bH^*\left(\bfX;(\Omega^{\ge p},\partial)\right)$. We note in passing 
that, just as in the case of $X$, the $LG$-action on $H^q(\Omega^p)$ is 
trivial (\cite{tel3}, Remark 8.10). 

\begin{maintheorem}\label{thickhodge}\hfill\begin{enumerate}
\item $H^\bullet(\bfX_\Sigma;\Omega^\bullet)$ is the free skew-commutative 
algebra generated by copies of $\Omega^0[\Sigma]$ and of $\Omega^1[\Sigma]$, 
in $H^m(\Omega^m)$, respectively $H^m(\Omega^{m+1})$, as $m$ ranges over the 
exponents of $\frg$.
\item The first Hodge-de Rham differential $\partial_1$ is induced by 
de Rham's operator $d:\Omega^0[\Sigma] \to\Omega^1[\Sigma]$ on generators, 
and spectral sequence collapses at $E_2$.
\end{enumerate}
\end{maintheorem}
\noindent The theorem will be proved in \S\ref{proof}. Assuming it, 
Proposition \ref{cohom} implies that $E_2$ already has the size of 
$H^\bullet(\bfX_\Sigma;\bC)$; this forces the vanishing of $\partial_2$ 
and higher differentials. 

\section{Uniform description of the cohomologies}\label{unified}
We now relate the Dolbeault and Macdonald cohomologies. In the process, we 
give a unified construction for the generating Dolbeault classes in Theorems 
\ref{thinhodge} and \ref{thickhodge}; this sets the stage for the proofs. 

\subsection{Moduli spaces and stacks.} In \S\ref{7}, we identified 
the thick flag variety $\bfX_\Sigma$ and the loop Grassmannian $X$ with the 
moduli spaces $\frM(\oSig,\wD)$ and $\frM(\bP^1,\bP^1\setminus \{0\})$ of 
$G$-bundles over the respective pairs. Their Dolbeault groups are described 
in Theorems \ref{thinhodge} and \ref{thickhodge}. For $\frM(\oSig)$, Hodge 
decomposition \cite{tel4} implies that $H^\bullet(\Omega^\bullet)$ is the 
free algebra on the bi-graded vector space $H^{\bullet,\bullet}(\oSig)^* 
\otimes\prim H^{\bullet,\bullet}(BG)$; this is Theorem \ref{cohom} with 
$D=\emptyset$. 

Consider the stack $\frM(\wD)$ of $G$-bundles on $\wD$. Such bundles are 
trivial (locally in any family), but their automorphisms are locally 
represented by the group $G[\wD]$ of regular formal loops. So $\frM(\wD)$  
is the classifying stack $BG[\wD]$. Cathelineau \cite{cat} identified the 
Hodge-de Rham sequence for the classifying stack of a complex Lie group 
$\cG$ (defined, say, from the simplicial realisation) with the holomorphic 
Bott-Shulman-Stasheff spectral sequence \cite{bs}
\begin{equation}\label{cat}
E_1^{p,q} = H^{q-p}(B\cG;\cO^{an}\otimes\sym^p \mathrm{Lie}(\cG)^*) 
			\Rightarrow H^{p+q}(B\cG;\bC),
\end{equation}
in which $E_1$ is the group cohomology with $\sym\mathrm{Lie}(\cG)^*
$-valued analytic co-chains and the abutment is the cohomology with 
constant coefficients. The result applies to any group sheaf $\cG$ over 
the site of algebraic or analytic spaces: indeed, \eqref{cat} is the 
\textit{descent spectral sequence} for the following fibration of 
classifying stacks, where $\widehat\cG$ denotes the formal group of $\cG$ at 
the identity:
\[
B\widehat\cG \hookrightarrow B\cG \twoheadrightarrow B(\cG/\widehat\cG).
\] 
The base of this fibration has the property that $H^\bullet(B(\cG/\hat\cG);
\cO)= H^\bullet(B\cG;\bC)$, by de Rham's theorem in the category of spaces. 
The first differential $H^n(\sym^p) \to H^{n-1}(\sym^{p+1})$ sends a group 
cocycle $\chi:\cG^{n+1}\to\sym^p$ to the sum of transposes of its derivatives 
$d_i\chi:\cG^n\times \mathrm{Lie}(\cG) \to\sym^p$ at $1$ along the 
components $i=0,\ldots,n$, symmetrised to land in $\sym^{p+1}$.

For simplicity, let $D$ be a single point, so $G[\wD]\cong G[[z]]$. 
Contractibility of $G[[z]]/G$ and the van Est sequence give an isomorphism 
\cite{tel3} between the cohomology $H^\bullet(BG[[z]];\sym^p\frg[[z]]^*)$ 
over the algebraic site of $BG[[z]]$ and the Lie algebra cohomology $H^
\bullet(\frg[[z]],\frg; \sym^p\frg[[z]]^*)$ (computed using continuous duals 
in the  Koszul complex; this is the restricted cohomology of \S\ref{2}). 
Theorem \ref{sym} then says that $E_1^{p,q}$ in the Hodge-de Rham sequence 
for $BG[[z]]$ is the algebra generated by the continuous duals of 
$\bC[[z]]dz$ and $\bC[[z]]dz$, in bi-degrees $(p,q)=(m,m)$ and $(m+1,m)$, 
respectively. The first differential converts an odd generator in $\Lambda
\otimes\sym$ to its even partner: this is induced by $s\mapsto 1$, or 
$n=0$ in \eqref{251}. We showed in \S\ref{2} that this leads to the dual 
of de Rham's operator, $d^*: (\bC[[z]]dz)^* \to\bC[z]^*$ on generators.    

\subsection{Sheaf cohomology for a pair.} \label{83}
For a coherent sheaf $\cS$ on $\oSig$, define the \textit{cohomology 
$H^\bullet(\oSig,\wD;\cS)$ relative to $\wD$} as the hyper-cohomology of 
the $2$-term complex $\cS\to\cS_{\wD}$, starting in degree $0$, mapping 
$\cS$ to its completion at $D$.\footnote{This is also the coherent sheaf cohomology \textit{with 
proper supports} on the open curve $\Sigma=\oSig \setminus D$; it only 
depends on the restriction of $\cS$ to $\Sigma$.} If $D =\emptyset$, this 
is the ordinary sheaf cohomology on $\oSig$; else, $H^0$ is the torsion 
of $\cS$ over $\Sigma$, and $H^1$ is identified with $\Hom_\oSig(\cS,
\Omega^1)^*$ by Serre duality. The groups relevant for us are  
\[\
H^1(\oSig,\wD;\cO) \cong \Omega^1[\Sigma]^*,\qquad
H^1(\oSig,\wD;\Omega^1) \cong \Omega^0[\Sigma]^*,
\] 
Serre dual to the opposite-degree differentials on $\Sigma$. Similarly, 
$H^\bullet(\oSig,\Sigma;\cS)$ is the hyper-cohomology of $\cS\to i_*i^*
\cS$, where $i:\Sigma\hookrightarrow\oSig$ is the inclusion. Again, we 
want $\cS = \Omega^{0,1}$, when $H^0$ vanishes and Serre duality describes 
the $H^1$'s as the continuous duals 
\[
H^1(\oSig,\Sigma;\cO) \cong \Omega^1[\wD]^*,
\qquad
H^1(\oSig,\Sigma;\Omega^1) \cong \Omega^0[\wD]^*,
\]
also known as the $\cO$- and $\Omega^1$-valued residues on $\oSig$ at $D$. 
When $\Sigma = \bP^1\setminus\{0\}$, these are the restricted duals of 
$\bC[z]dz$ and $\bC[z]$. 

The following summarises Theorems \ref{sym}, \ref{thinhodge} and \ref
{thickhodge}. 

\begin{theorem}\label{84}
Let $S$ stand for $\wD$, $\oSig$ or one of the pairs $(\oSig,\Sigma)$ 
or $(\oSig,\wD)$. Then, $\bH^q(\frM(S);\Omega^p)$ is the free skew-commutative 
algebra on 
\[
H^\bullet(S;\Omega^\bullet)^*\otimes \gen^{\bullet,\bullet}(BG).
\] 
The first Hodge-de Rham differential $\partial_1$ is induced by de Rham's 
differential on generators, and all higher differentials vanish. \qed
\end{theorem} 

\subsection{Dolbeault generators.} \label{85}
We are now in a position to describe the generating Dolbeault classes. 
For a principal $G$-bundle $\cP$ over a base $B$, the tangent bundle to 
the total space of $\cP$ is $G$-equivariant and descends thus to the base, 
where it gives an extension $\ad_\cP \to T\cP/G \to TB$. This defines the 
\textit{Atiyah class} in $H^1(B; \ad_\cP\otimes \Omega^1)$. With $S$ as in 
\eqref{84} and the universal $G$-bundle $\cP$ over $S\times\frM(S)$, we 
obtain the universal Atiyah class  
\[
\alpha_S\in \bH^1\left(S\times\frM(S);\ad_\cP\otimes\Omega^1\right),
\]
keeping in mind that differentials form a complex when $\frM$ is a stack. 
An invariant polynomial $\Phi$ of degree $d$ on $\frg$ defines a class 
$\Phi(\alpha)\in \bH^d\left(S\times\frM(S);\Omega^d\right)$. 

\begin{proposition}\label{86}
Serre duality contraction of $\Phi(\alpha)$ with $H^a(S;\Omega^b)^*$ gives 
the Dolbeault generators in $\bH^{d-a}\left(\frM(S);\Omega^{d-b}\right)$ of 
Theorem \ref{84}, as $\Phi$ ranges over $\gen^{d,d}(BG)$.
\end{proposition} 

\begin{proof}
Clearly, the construction is a Dolbeault refinement of the topological 
transgression in Theorem \ref{cohom}. For $S=\oSig$, there is nothing left 
to show, since Hodge decomposition on $\oSig$ and $\frM$ equates Dolbeault 
and de Rham cohomologies. For $S= (\oSig,\Sigma)$ or $(\oSig,\wD)$ (the 
flag varieties), the proposition completes the statement of Theorems \ref
{thinhodge} and \ref{thickhodge}, and will be proved in the next section. 
For the remaining case $S=\wD$, we now relate the newly constructed 
generators to those of Theorem \ref{sym}. For simplicity we let $\wD=
\mathrm{Spf}\,\bC[[z]]$.

The cotangent complex of $\wD\times BG[\wD]$ splits as $\bC dz\oplus
\frg[[z]]^*[-1]$, the second summand being a bundle over $BG[[z]]$ under 
the co-adjoint action. We then have
\begin{equation}\label{87}
H^1\left(\wD\times BG[[z]];\ad\otimes\Omega^1\right) = 
	H^1\left(BG[[z]];\frg[[z]]dz\right)\oplus
		\Hom_{G[[z]]}\left(\frg[[z]];\frg[[z]]\right).
\end{equation}
We shall show at the end of this section that the two components of $\alpha$ 
are the group co-cycle $\gamma\mapsto -d\gamma\cdot\gamma^{-1}$ and the 
identity map $\mathrm{Id}$. (Both groups are in fact free $\bC[[z]]$-modules 
of rank one, generated by the respective classes, but we do not need this 
fact.) Applying $\Phi$ to the second component of $\alpha$ leads to  
\[
\Phi(\mathrm{Id})\in 
	\Hom_{G[[z]]}\left(\hat\sym{}^d\frg[[z]];\bC[[z]]\right)
\]
which is just the point-wise application of $\Phi$. Contracting with a 
Fourier mode $z^n\in\bC[[z]]$ gives the co-cycle $S(-n)$ of Theorem 
\ref{sym}, viewed as an element of 
\[
H^d(BG[[z]]; \Omega^d) = \left(\sym^d\frg[[z]]^*\right)^{G[[z]]}.
\]
The first factor in \eqref{87} squares to $0$, because $\Omega^1(\wD)$ 
does so; hence, $\Phi$ takes no more than one entry from there. 
Absorbing one entry from $T^*\wD\otimes\ad = \frg[[z]]dz$, followed by 
contraction against $z^{n-1}dz$, leads to a group $1$-cocycle $G[[z]]\to 
\hat\sym^{d-1}\frg[[z]]^*$---the contraction of $\Phi$ with $-d\gamma 
\cdot\gamma^{-1}$. Via the van Est isomorphism with $H^1(\frg[[z]],
\frg;\;.\;)$, this becomes the odd generator $E(-n)$.
\end{proof}

\subsection{The Atiyah class.}\label{88}
To complete the proof, we must say more about $\alpha$. Decompose $\Omega^1 
=\Omega^1_S\oplus\Omega^1_{\,\frM(S)}$; the two components of $\alpha$ can 
be interpreted as the Kodaira-Spencer deformation maps for $\cP$, first 
regarded as a family of bundles over $\frM(S)$ parametrised by $S$, and 
then as a family of bundles on $S$ parametrised by $\frM(S)$. Now, $\Omega
^1_{\,\frM(S)}$ is the complex $R\pi_*\ad_\cP[1]$ for the projection $\pi$ 
along $S$ to $\frM(S)$, and from the definition of $\frM(S)$, $\alpha$ 
has a tautological component 
\[
\mathrm{Id}\in R\Hom_{\frM(S)}
		\left(R\pi_*(\ad_\cP);R\pi_*(\ad_\cP)\right)\cong
	\bH^1\left(S\times\frM(S); \Omega^1_{\,\frM(S)}\right).
\]
The more geometric component is in $H^1(\Omega^1_S\otimes \ad_\cP)$; its 
leading term in $\Gamma\left(S;\Omega^1\otimes R^1\pi'_*\ad_\cP\right)$, 
where $\pi'$ is the projection to $S$, represents the local Kodaira-Spencer 
deformation map for $\cP$. 

\begin{remark} In our examples, the remaining information in $H^1\left
(S;\Omega^1\otimes\pi'_*\ad_\cP\right)$, is nil: for $S=\oSig, (\oSig,
\Sigma)$ or $(\oSig,\wD)$, the $\pi'_*$-sheaf is null, whereas if $S=\wD$, 
$H^1=0$.
\end{remark}

Let us spell out $\alpha_{\wD}$. $T\cP/G$ is the complex over $\wD\times 
BG[\wD]$ descended from $T(G\times\wD)$, with the $G[\wD]$-translation action. 
Postponing for a moment the matter of the $G[\wD]$-action, the underlying 
complex is $\frg[\wD]\to (\frg\oplus T\wD)$; the arrow, coming from the 
fibre-wise translation action, is the evaluation map $\frg[\wD]\times\wD\to 
\frg$. This is the tautological component of $\alpha$, becoming the identity 
in the second summand in \eqref{87}. Now, $G[\wD]$ acts by ad on the first 
term of the complex, whereas the second term is an extension of 
$G[\wD]$-equivariant bundles
\[
\left[\frg\to \left(T\cP|_{\wD}\right)/G \to T\wD \right]
	\in \Ext^1_{\wD\times BG[\wD]}(T\wD;\frg).
\] 
A $\gamma\in G[\wD]$ changes a splitting $\frg\oplus T\wD$ by sending a 
section $(\xi,v)$ to $(\gamma\xi\gamma^{-1}-v(\gamma)\cdot\gamma^{-1},v)$, 
and the derivative term represents the class of $\gamma\mapsto -d\gamma 
\cdot\gamma^{-1}$ in the $\Ext$ group.

\section{Proof of Theorems \ref{thinhodge} and \ref{thickhodge}}
\label{proof}
We now compute the Dolbeault cohomology for thick flag varieties. For 
convenience, in this section we write $\bfX$ for $\bfX_\Sigma$ and $\frM$ 
for $\frM(\oSig)$, and continue to assume that $D$ is a single point; the 
changes needed for the general case are obvious. A small modification then 
gives us Theorem~\ref{thinhodge}.  

\subsection{Setting up the spectral sequence.} \label{91}
Uniformisation (\S\ref{71}) realises $\frM$ as the quotient stack 
$G[[z]]\setminus \bfX$. Equivariance under the translation $G[[z]]$-action 
on $\bfX$ makes the bundle $\Omega^p$ of differential $p$-forms descend 
to a bundle on $\frM$; we denote the descended bundle by $\Omega^p_\bfX$. 
The complex of differentials $\Omega^r=\Omega^r_\frM$ on $\frM$ is 
represented by a Koszul-style complex of bundles
\begin{equation}\label{92}
\Omega^r \sim \left(\Omega^r_\bfX\xrightarrow{\kappa} 
		\sym^1\frg[[z]]^*\otimes\Omega^{r-1}_\bfX\xrightarrow{\kappa}	
		\sym^2\frg[[z]]^*\otimes\Omega^{r-2}_\bfX\xrightarrow{\kappa}
		\cdots\right),
\end{equation}
cohomologically graded by symmetric degree. To describe the differential, 
observe that a choice of $\gamma\in LG$ identifies the tangent space to
$\bfX$ at $\gamma G[\Sigma]$ with $L\frg/\frg[\Sigma]$; thereunder, 
$\kappa$ at $[\gamma]= G[[z]]\gamma G[\Sigma]\in \frM$ is induced by 
the $\gamma$-twisted dual to the natural projection $\frg[[z]]\to L\frg/
\frg[\Sigma]$.

The complex \eqref{92} has finite length, so it leads to a convergent 
spectral sequence with 
\begin{equation}\label{93}
E_1^{k,l} = H^l\left(\frM; \sym^k\frg[[z]]^*\otimes \Omega^{r-k}
		_\bfX\right) \Rightarrow H^{k+l}\left(\frM;\Omega^r\right).
\end{equation}
There is one such spectral sequence for each $r\ge 0$, but the product, 
which is compatible with the differentials, mixes them. We have an 
identification of cohomologies 
\[
H^l(\frM;\sym^k\otimes\Omega_\bfX^{r-k}) = 
		H^l_{G[[z]]}\left(\bfX; \sym^k\otimes\Omega^{r-k}\right),
\] 
where $H^k_{G[[z]]}$ is the (algebraic) equivariant cohomology. 

\subsection{The Key Factorisation.} \label{key}
Our $E_1$ term \eqref{93} factors as 
\begin{equation}\label{94}
E_1^{k,l} = \bigoplus\nolimits_s H^s_{G[[z]]}\left( 
			\sym^k\frg[[z]]^*\right)
			\otimes H^{l-s}\left(\bfX; \Omega^{r-k}\right).
\end{equation}
\textit{A priori}, the right-hand side is the ${}^LE_2^{s,l-s}$ term in 
the Leray sequence for the sheaf $\sym^k\otimes\Omega^{r-k}$ and the 
morphism $\frM\to BG[[z]]$. However, no differentials are present, 
because ${}^LE_2$ is generated from the bottom edge ${}^LE_2^{s,0}$ by 
cup-product with classes which live on the total space: indeed, because 
$G[[z]]$ acts trivially on the cohomology and $H^{>0}_{BG[[z]]}(\cO)=0$, 
we have an isomorphism
\[
H^l\left(\bfX; \Omega^{r-k}\right) \cong 
		H_{G[[z]]}^l\left(\bfX; \Omega^{r-k}\right).
\]
This also shows that \eqref{94} is a natural isomorphism, and not just 
the $\Gr$ of one.
 
\subsection{Determining the spectral sequence.}
The factor $H^s(BG[[z]]; \sym^k \frg[[z]]^*)$ is isomorphic to the 
Macdonald cohomology of Theorem B. The abutment $H^s\left(\frM; \Omega^r
\right) = H^{r,s}(\frM;\bC)$ is also known \eqref{cohom}. We now construct an 
obvious candidate for the spectral sequence, with a map to \eqref{93}, and 
prove by induction on $r$ that the obvious candidate is correct. This last 
argument is a variation on Zeeman's comparison theorem \cite{z}. 

\begin{proposition}\label{95}
The sum over all $r$ of the spectral sequences \eqref{93} is the commutative 
differential bi-graded algebra freely generated by copies of the 
differential bi-graded vector spaces $\Omega^0[\Sigma] \to\Omega^1[\wD]^*$, 
in bi-degrees $(k,l)=(0,m)$ and $(m,1)$, and $\Omega^1[\Sigma]\to \Omega^0 
[\wD]^*$ in bi-degrees $(k,l)=(0,m)$ and $(m+1,0)$, respectively, as $m$ 
ranges over the exponents of $\frg$.

The arrows above are dual to the connecting maps $H^0(\wD;\Omega^i)\to 
H^1(\oSig,\wD;\Omega^i)$, $i=1,0$. 
\end{proposition}
\noindent Concretely,  
\[
\Omega^1[\wD]^*= \bC((z))/\bC[[z]],\qquad \Omega^0[\wD]^* = \bC((z))/\bC[[z]]\otimes dz, 
\]
and the maps are the principal parts at $z=0$ on $\oSig$ (cf.~\S\ref{83}). 
Here is the location of the generators, with respect to the decomposition 
\eqref{94}:
\[
\begin{tabular}{lccccc}
 space &\vline & $k$&$l$&$r$&$s$ \\
\hline
$\Omega^0[\Sigma]$  &\vline& $0$ & $m$ & $m$ & $0$ \\ 
$\Omega^1[\Sigma]$  &\vline& $0$ & $m$ & $m+1$ & $0$\\ 
$\Omega^1[\wD]^*$  &\vline & $m$ & $1$ & $m$ & $1$ \\ 
$\Omega^0[\wD]^*$  &\vline & $m+1$ & $0$ & $m+1$ & $0$ \\
\end{tabular}
\]
The spectral sequence differential which originates at $\Omega^i[\Sigma]$ 
has length $m+i$. 

\begin{proof}
The candidate generators are mapped to $E_1$ as in Proposition \ref{86}. 
We will see in \S\ref{98} below that the terms $\Omega^i [\Sigma]$ 
survive to $E_{m+i-1}$, and that the differential $\delta_{m+i}$ maps them 
into the $\Omega^{1-i}[\wD]^*$ in the way indicated. This good behaviour 
is enforced by the tautological summand in the Atiyah class in \S\ref{88}. 
Observe also that the kernels and co-kernels of these differentials are 
the Dolbeault groups of $\oSig$, so the fact that they define the generating 
classes for $H^s\left(\frM;\Omega^r\right)$ (and therefore survive to 
$E_\infty$) is already known from the Hodge decomposition of $\frM$ \cite
{tel4}, and the topological (Atiyah-Bott) construction of its cohomology generators.   

Let now ${}'E_{n}^{k,l}$, $n\ge 1$,  be the spectral sequence with 
multiplicative generators and differentials as in \eqref{95}. We will 
show by induction on $r$ that the map to $E_n^{k,l}$ we constructed 
is an isomorphism. For $r=0$, this merely says that $H^0(\bfX;\cO) = 
\bC$ and $H^{>0}(\bfX;\cO)=0$, which was shown in \cite{tel3}. If 
the assumption holds up to $r$, then the multiplicative decomposition 
\eqref{94} shows that, for $r+1$, ${}'E_1^{k,l}\cong E_1^{k,l}$, except 
perhaps on the left edge $k=0$. 

The assumption also implies that the spectral sub-sequence of ${}'E_n^
{k,l}$, $k>0$, obtained by deleting the left edge, converges to the 
hyper-cohomology of the sub-complex of $\Omega^{r+1}$ 
\[
\Omega^{r+1}_+:= \sym^1\frg[[z]]^*\otimes \Omega^r_\bfX 
	\xrightarrow{\kappa}\sym^2\frg[[z]]^*\otimes \Omega^{r-1}_\bfX
	\xrightarrow{\kappa}\cdots.
\]
Our construction gives a map between the long exact sequences of 
cohomologies over $\frM$,
\[
 \ldots\to\bH^l\left(\frM; \Omega^{r+1}_+\right) 
 	\to \bH^l\left(\frM;\Omega^{r+1}\right) \to {}'E_1^{0,l} 
	\to \bH^{l+1}\to\ldots,
\]
obtained from the spectral sub-sequence, and 
\[
 \ldots\to\bH^l\left(\frM; \Omega^{r+1}_+\right) 
 	\to \bH^l\left(\frM;\Omega^{r+1}\right) 
	\to H^l_{G[[z]]}\left(\bfX;\Omega^{r+1}\right)
	\to \bH^{l+1}\to\ldots
\]
arising from the sub-complex. As explained in \eqref{94}, we can omit 
the $G[[z]]$-subscript in the third cohomology, and the Five Lemma gives 
the desired isomorphism ${}'E_1^{0,l}\cong H^l\left(\bfX;\Omega^{r+1}
\right)$.
\end{proof}

\subsection{The Hodge differentials.}
Since the construction of generators is compatible with de Rham's operator, 
the first Hodge-de Rham differentials are as described in 
Theorem~\ref{thickhodge}. 

\subsection{The loop Grassmannian.}\label{97}
To prove Theorem \ref{thinhodge}, we repeat the argument above, but use
the presentation $\frM(\bP^1) = G[z^{-1}]\setminus X$ of the stack of 
$G$-bundles. The complex \eqref{92} representing the differentials is now 
a pro-vector bundle, completed for the $z^{-1}$-adic filtration on 
$\sym^k\frg[z^{-1}]\otimes\Omega_X^{r-k}$. The key factorisation result 
\eqref{94} continues to apply (completed in the filtration topology), this 
time using Proposition \ref{62}. 

\subsection{Leading differentials in $E_\bullet$.} \label{98}
We now check the good behaviour of the Dolbeault generators assumed in 
the proof. The argument is a convoluted tautology, but we include it 
nonetheless for completeness. An invariant $\Phi\in\sym^{m+1}\frg^*$, 
applied to the Atiyah class
\[
\alpha_\frM\in \bH^1\left(\oSig\times\frM; 
		\Omega^1_\oSig\otimes\ad_\cP\oplus
		(\Omega^1_\bfX\to\sym^1\frg[[z]]^*)\otimes\ad_\cP\right),
\]
accepts, for dimensional reasons, a single non-tautological (first) entry. 
For the same reason, this entry will be detected under contraction with 
the first set of generators $\Omega^0[\Sigma]$ in Proposition \ref{95}, 
but will be killed by the second. So the first family of generators contain 
the tautological summand to degree $m$; the second, to degree $m+1$. 

We project the tautological (second) component of $\alpha_\frM$ to $\bH^1\left(\oSig\times\frM; \Omega^1_\bfX\otimes\ad_\cP\right)$. 
Lifted to $\bfX$, this is the tautological component of $\alpha_\bfX$, 
and the two classes have the common refinement 
\begin{equation}\label{99}
\mathrm{Id}\in R\Hom_\frM\left(R\pi_*\ad_\cP; R\pi_*\ad_\cP\right)
	\cong\bH^1\left((\oSig,\wD)\times\frM; 
			\Omega^1_\bfX\otimes\ad_\cP\right)
\end{equation}
(notations as in \S\ref{88}). Note that cup-product of \eqref{99} with 
classes in $H^\bullet(\oSig\times\ldots)$ lands in $H^\bullet((\oSig,\wD)
\times\ldots)$, and such classes can be contracted with (= integrated 
against) functions and forms on $\Sigma$. \textit{Consequently, contraction 
of $\Phi(\alpha_\Sigma)$ with $\Omega^i[\Sigma]$ gives well-defined classes 
in the truncated complex $\bH^m(\frM; \Omega^{m+i}/\sym^{m+i})$}. In 
particular, the Dolbeault generators in Proposition \ref{95} survive to 
$E_{m+i-1}$. 

To conclude, we must identify the differentials $\delta_{m+i}$. These arise 
from the failure of 
\[
\Phi(\alpha)\in\bH^{m+1}\left((\oSig,\wD)\times\frM; \Omega^{1-i}_\oSig
		\otimes\Omega^{m+i}/\sym^{m+i}\right)
\] 
to lift to $\bH^{m+1}(\,.\,; \Omega^{1-i}_\oSig\otimes\Omega^{m+i})$. The 
obstruction is detected by a connecting homomorphism to $H^{2-i}(\,.\,; 
\Omega^{1-i}_\oSig \otimes\sym^{m+i})$ (of degree $1$, but we have shifted 
degrees because of $\sym$). Contraction leads to our differentials, which 
land in
\begin{equation}\label{9-10}
H^{1-i}\left(\frM; \sym^{m+i}\frg[[z]]^*\right) \cong 
		H^{1-i}\left(BG[[z]]; \sym^{m+i}\frg[[z]]^*\right);
\end{equation}
we used the key factorisation \eqref{key} for the isomorphism. 

We can identify the connecting map in a different way. The class $\Phi
(\alpha)$ \textit{does} lift to the full differentials $\Omega^{m+i}$, 
but only over $\oSig\times\frM$. A diagram chase then shows that our 
obstruction is the image of the restricted $\Phi(\alpha)$ under the 
connecting map 
\[
H^{1-i}\left(\wD\times\frM; \Omega^{1-i}_\oSig \otimes\sym^{m+i}
			\frg[[z]]^*\right) 
	\xrightarrow{\:\partial\:}
H^{2-i}\left((\oSig,\wD)\times\frM; \Omega^{1-i}_\oSig \otimes\sym^{m+i}
			\frg[[z]]^*\right).
\]
However, the universal bundle $\cP$ over $\oSig\times\frM$, when restricted 
to $\wD\times\frM$, is pulled back from the universal bundle on $\wD\times
BG[\wD]$; hence, so is its Atiyah class and $\Phi(\alpha)$, and we can 
replace $\frM$ with $BG[\wD]$ above. Moreover, the $\partial$'s are the 
residue maps appearing in Proposition \ref{95}. This identifies the $\delta
_{m+i}$ with the Dolbeault classes describe in Proposition \ref{95}. 

\section{Related Lie algebra results} \label{10}
\subsection{Dolbeault cohomology as Lie algebra cohomology.} We now give 
a Lie algebra interpretation of the Dolbeault cohomology of 
the loop Grassmannian $X= LG/G[[z]]$. The dual of $\frg((z))/\frg[[z]]$ 
is identified with $\frg[[z]]dz$ by the residue pairing. The $p$-forms on 
$X$ are then sections of the pro-vector bundle associated to the adjoint 
action of $G[[z]]$ on $\hat\Lambda^p\frg[[z]]dz$. Recall that the 
latter is $z$-adic completion of the exterior power. For modules thus 
completed, it is sensible to form the \textit{continuous} $\frg[[z]]
$-cohomology, resolved by the Koszul complex of \textit{continuous} linear 
maps 
\begin{equation}\label{101}
\mathrm{Hom}\left(\hat\Lambda^\bullet\frg[[z]];\hat\Lambda^p
(\frg[[z]]dz)\right);
\end{equation}
in this case, we get the \textit{inverse limit} of cohomologies\footnote
{Finite-dimensionality of cohomology shows that there are no $R^1\lim$ 
terms to worry about.} of the $z$-adic truncations of the coefficients. 
We emphasise, however, that the complex \eqref{101} has infinite-dimensional 
$(z,\frg)$-eigenspaces, which is a serious obstacle to a direct computation 
of its cohomology as in Chapter I. 

\begin{proposition}
The Lie continuous algebra cohomology $H^q_{cts}\left(\frg[[z]],\frg;
\hat\Lambda^p\frg[[z]]dz\right)$ resolved by (\ref{101}) is naturally 
isomorphic to $H^q(X;\Omega^p)$. 
\end{proposition}
\begin{proof}
Contractibility of $G[[z]]/G$ gives a natural ``van Est" isomorphism 
\cite{tel3} 
\[
H^q_{cts}\left(\frg[[z]],\frg;\hat\Lambda^p\frg[[z]]dz\right) =
H^q_{G[[z]]}\left(\hat\Lambda^p\frg[[z]]dz\right);
\]
noting that the $H^q\left(X;\Omega^p\right)$ are the $q$th derived 
functors of induction from $G[[z]]$ to $LG$, the group and Dolbeault 
cohomologies are related by Shapiro's spectral sequence
\[
E_2^{r,s} = H^r_{LG}\left(H^s\left(X;\Omega^p\right)\right)
\Rightarrow H^{r+s}_{G[[z]]}\left(\hat\Lambda^p(\frg[[z]]dz)\right).
\]
Alternatively, this is the Leray sequence for the morphism $BG[[z]]\to 
BLG$, with fibre $X$. Either way, $H^s(X;\Omega^p)$ is a trivial $LG$-module, 
so its higher $LG$-cohomology vanishes; the spectral sequence collapses 
and we obtain the asserted equality.
\end{proof}

\subsection{Thick flag varieties.} An obvious variation replaces the 
formal disk $\mathrm{Spf}\bC[[z]]$ by a smooth affine curve $\Sigma$. 
We consider the Lie algebra cohomology $H^q\left(\frg[\Sigma]; 
\Lambda^p\,\Omega^1(\Sigma;\frg)\right)$. The answer carries now a 
contribution from the non-trivial topology of the group $G[\Sigma]$. 
As in \cite{tel3}, the van Est sequence collapses at $E_2$, leading, 
by the same argument as above, to
\begin{proposition}
$H^\bullet\left(\frg[\Sigma]; \Lambda^p\,\Omega^1
(\Sigma;\frg)\right) \cong H^\bullet(\bfX_\Sigma;\Omega^p)\otimes 
H^\bullet(G[\Sigma];\bC)$, naturally.\qed
\end{proposition}
\noindent The homotopy equivalence of $G[\Sigma]$ with the corresponding 
group of continuous maps shows that the topological factor $H^*(G[\Sigma]; 
\bC)$ is isomorphic to 
\begin{equation}\label{10-6}
H^\bullet(G;\bC)\otimes H^\bullet(\Omega G;\bC)^{1-\chi(\Sigma)},
\end{equation}
with $\Omega G$ denoting the space of based continuous loops and $\chi$ the
Euler characteristic. This is also isomorphic to the Lie algebra cohomology $H^\bullet(\frg[\Sigma];\bC)$ \cite{tel3}.

\subsection{Strong Macdonald for smooth curves.} The method of \S\ref
{proof} allows us to carry out the long-postponed proof of the 
higher-genus version of the strong Macdonald theorem. 
\begin{proof}[Proof of Theorem \ref{1-15}.]
We use the construction of \S\ref{proof}, but realise the moduli stack 
$\frM$ of $G$-bundles on $\oSig$ as the quotient $X/G[\Sigma]$, and 
present the differentials $\Omega^r_\frM$ by a complex of pro-vector 
bundles 
\[
\Omega_X^{r-k}\widehat\otimes\,\hat\sym{}^k\frg[\Sigma]^*,
\] 
The factorisation replacing \eqref{94} now reads 
\[
E_1^{k,l} = \bigoplus\nolimits_s H^s_{G[\Sigma]}\left( 
			\hat\sym{}^k\frg[\Sigma]^*\right)
			\widehat\otimes H^{l-s}\left(X; \Omega^{r-k}\right);
\]
the Dolbeault cohomologies of $X$ being known, the desired group 
cohomologies are again determined by induction on $r$, with the difference 
that it is the right $(r+1)$st edge of the sequence that is new, in the 
inductive step. Collapse of the van Est sequence leads to the factor 
$H^\bullet(\frg[\Sigma];\bC)$ when switching from group to Lie algebra
cohomology.
\end{proof}

\part{Positive level}

The loop group $LG$ admits central extensions by the circle; when $G$ is 
semi-simple, these are parametrised up to isomorphism by a \textit{level} 
in $H^3_G(G;\bZ)$.\footnote{There are additional choices for the torus 
factors, but only one of them is interesting \cite{pseg}.} When $G$ is 
simple and simply connected, $H^3_G(G;\bZ)\cong \bZ$, and positive levels 
lead to the interesting class of \textit{highest-weight representations} 
of $LG$, also called \textit{integrable highest-weight modules} of $\frg
((z))$. These have a Borel-Weil realisation as spaces of sections of vector 
bundles over the genus-zero thick flag variety $\bfX$ (\S\ref{11-7}), and 
carry a semi-simple $\bC^\times$-action intertwining with the $z$-scaling. 
The eigenvalues of its infinitesimal generator, called \textit{energies} or 
\textit{$z$-weights}, are bounded above.   

In \S\ref{11} gives a positive-level analogue of Theorem \ref{sym}, 
which includes in the coefficients of Macdonald cohomology a highest-weight 
$LG$-representation $\bfH$. This entails the \textit{vanishing of higher 
cohomology}. As a result, the analogue of Macdonald's constant term---the 
$(z,s)$-weighted Lie algebra Euler characteristic in \eqref{14}---refines 
the $z$-dimension of the $G$-invariant part of $\bfH$, detecting an affine 
analogue of \textit{R.~Brylinski's filtration} \cite{bryl}, originally 
defined on weight spaces of $G$-representations (Remark \ref{11-0}).

Central extensions of $LG$ lead to algebraic line bundles over the 
loop Grassmannians $X$ and $\bfX$. The sections of the level $h$ line 
bundle $\cO(h)$ over $\bfX$ span the highest-weight \textit{vacuum 
representation} $\bfH_0$. In \S\ref{twist}, we determine the level 1 
Dolbeault cohomologies $H^q(\bfX;\Omega^p(1))$, for simply laced $G$. A 
combinatorial application is given in \S\ref{bailey}. 

\section{Brylinski filtration on loop group representations}\label{11}
Let $\bfH$ be a highest-weight representation of $LG$; it is the 
\textit{direct sum} of its $z$-weight spaces $\bfH(n)$. We assume that 
the level is positive on each simple or central factor of $\frg$; the 
only level-zero representation is trivial, and has been discussed 
already.

\begin{maintheorem}\label{vanish}
$H^k(\frg[[z]],\frg; \bfH \otimes \sym^p\frg[[z]]^*)$ vanishes for 
positive $k$.
\end{maintheorem}
\noindent With respect to Chapter I, this is the restricted cohomology 
$H^k_{res}(\frg[z],\frg; \bfH \otimes \sym^p\frg[z]^*_{res})$.
\begin{proof} 
For abelian $\frg$, $\bfH$ is a sum of Fock representations, and so is 
injective for $(\frg[z],\frg)$. Assume now that $\frg$ is simple and $\bfH$ 
has level $h>0$. In the notation of \S\ref{3}, with the operator $\bar
\partial$ on $\bfH\otimes\Lambda\otimes\sym$ modified to include the 
$\frg[z]$-action on $\bfH$, Theorem 2.4.7 from \cite{tel1} becomes 
\[
\overline\square = \square + T_S^\Lambda + k\cdot(1+\frac{h}{2c}) =
	 \square + D + k\cdot h/2c, 
\]
where the second identity follows as in Proposition \ref{317}. This is 
strictly positive if $k>0$. 
\end{proof}

\begin{remark} \label{11-0}
Theorem \ref{vanish} has finite-dimensional analogues for $G$-modules $V$
and the Borel and Cartan sub-algebras $\frb,\frh\subset\frg$: the 
higher cohomologies $H^{>0}(\frb,\frh; V\otimes\sym\frb^*)$ and $H^{>0}
(\frb,\frh; V\otimes\sym\frn^*)$ vanish ($\frn=[\frb,\frb]$). Using the 
Peter-Weyl decomposition of the polynomial functions on $G$, this is 
equivalent to the vanishing of higher cohomology of $\cO$ over $G\times_B 
\frb$ and $T^*G/B= G\times_B \frn$, and follows from the 
Grauert-Riemenschneider theorem (cf.\ the proof of Lemma \ref{412}).
\end{remark}

\subsection{Shift in the grading.} \label{11-1}
For reasons that will be clear below, we now replace $\frg[[z]]$ in the 
symmetric algebra by the differentials $\frg[[z]]dz$. This does not alter 
the $\frg[[z]]$-module structure, but shifts $z$-weights by the symmetric 
degree. To match the usual conventions, we set $q=z^{-1}$ and consider the 
$q$-Euler characteristic in the restricted Koszul complex \eqref{34}, 
capturing the symmetric degree by means of a dummy variable $t$. After our 
shift, the isomorphism in Theorem \ref{sym} leads to the following identity, 
where $\mathrm{CT}$ denotes the $G$-constant term, after expanding the 
product into a formal $(q,t)$-series with coefficients in the 
representation ring of $G$:
\begin{equation}\label{11-2}
\mathrm{CT}\left[\bigotimes_{n>0}\frac{1-q^n\cdot\frg}
			{1-tq^n\cdot\frg}\right]
		= \prod_{k=1}^\ell\prod_{n > m_k}
			\frac{1-t^{m_k} q^n}{1-t^{m_k+1}q^n},
\end{equation}

\subsection{Constant term at positive level.} The \textit{$q$-dimension} 
$\dim_q \bfH := \sum\nolimits_n q^{-n}\dim\bfH(n)$, convergent for $|q|<1$, 
captures the $z^{-1}$-weights. Using the Koszul resolution of cohomology, 
Theorem \ref{vanish} equates the $q$-dimensions of the invariants with a 
$G$-constant term,
\begin{equation}\label{11-4}
\mathrm{CT}\left[\bfH\otimes\bigotimes_{n>0}
		\frac{1-q^n\frg}{1-tq^n\frg}\right] =
\sum_{p\ge 0} t^p\dim_q
	\left\{\bfH\otimes\sym^p(\frg[[z]]dz)^*\right\}^{\frg[z]}.
\end{equation}
When $G$ is simple and $\bfH$ is irreducible, with highest energy zero and 
highest weight $\lambda$, the Kac character formula \cite{kac} converts 
the $q$-representation $\bfH\otimes\prod_{n>0}(1-q^n\frg)$ of $G$ into the 
sum  
\begin{equation}\label{11-5}
\sum_{\mu\in \lambda +(h+c)P} 
		\pm q^{\frac{c(\mu)-c(\lambda)}{h+c}}V_\mu,
\end{equation}
where $c(\mu) = (\mu+\rho)^2/2$, $c$ is the dual Coxeter number of $\frg$, 
$P\subset \frh^*$ the integer lattice and $\pm V_\mu$ is the signed 
$G$-module induced from the weight $\mu$ (the sign depending, as usual, 
on the Weyl chamber of $\mu+\rho$). So the left side in \eqref{11-4} is 
also \begin{equation}\label{11-6}
\sum_{\mu\in \lambda+ (h+c)P} \pm q^{\frac{c(\mu)-c(\lambda)}{h+c}}
	\,\mathrm{CT}\left[\frac{V_\mu}{
		\bigotimes\nolimits_{n>0}(1-tq^n\frg)}\right].
\end{equation}
Its analogy with Macdonald's constant term becomes compelling, if we use 
the Kac formula at $h=0$ to equate the left side in \eqref{11-2} with the 
sum \eqref{11-6} for $\lambda=h=0$.

\subsection{Brylinski filtration.} \label{11-7}
Recall the Borel-Weil construction of $\bfH$. The \textit{thick flag variety} 
(\S\ref{7}) $\bfX:= G((z))/G[z^{-1}]$ of the formal Laurent loop group 
$G((z))$ carries the level $h$ line bundle $\cO(h)$ and the vector bundle 
$\cV_\lambda$, the latter defined from the action of $G[z^{-1}]$ on $V_
\lambda$ by evaluation at $z=\infty$. Then, $\bfH$ is the space of algebraic 
sections of $\cV_\lambda(h):=\cV_\lambda\otimes\cO(h)$ over $\bfX$. 

Restricted to the \textit{big cell} $\bfU\subset\bfX$, the orbit of the 
base-point under $G[[z]]$, $\cV_\lambda(h)$ is trivialised by the action of 
the subgroup $\exp(z\frg[[z]])$. Now, $\bfU$ is identified with $\frg[[z]]dz$ 
by $G[[z]]/G\ni\gamma \mapsto d\gamma\cdot \gamma^{-1}$, and the resulting 
affine space structure is preserved by the left translation action of $G[[z]]$. Sections of $\cV_\lambda(h)$, having been identified with $V_\lambda$-valued polynomials, are increasingly filtered by their degree, and this gives an increasing, $G[[z]]$-stable filtration of $\bfH$.
\begin{theorem}\label{11-8}
We have a natural isomorphism $\mathrm{Gr}_p\,\bfH^G \simeq 
	\left\{\bfH\otimes\sym^p(\frg[[z]]dz)^*\right\}^{\frg[z]}$.
\end{theorem}
\begin{proof}
With the conjugation action of $G[[z]]$, $\sym(\frg[[z]]dz)^*$ is the 
associated graded space of $\bC[\bfU]$, the space of polynomials on the 
open cell $\bfU\simeq\frg[[z]]dz$, filtered by degree. In the Borel-Weil 
realisation, $\bfH\otimes\bC[\bfU]$ is a subspace of the $V$-valued 
functions on $\bfU\times\bfU$, filtered by the degree on the second factor. 
It follows that the subspace of invariants under the diagonal translation 
action of $\frg[z]$ gets identified, by restriction to the first $\bfU$, 
with the $G$-invariants in $\bfH$, endowed with the Brylinski filtration. 
Cohomology vanishing gives an isomorphism
\[
\mathrm{Gr}_p\left\{\bfH\otimes\bC[\bfU]\right\}^{\frg[z]} = 
\left\{\bfH\otimes\sym^p(\frg[[z]]dz)^*\right\}^{\frg[z]},
\] 
leading to our theorem.
\end{proof}

\begin{remark}\label{11-10}
Applied to a $G$-representation $V$ and the cohomology vanishing in 
Remark \ref{11-0}, the same argument defines Brylinski's filtration on 
the zero-weight space $V^\frh\cong (V\otimes\sym\frn^*)^\frb$.
\end{remark}

\subsection{The basic representation.} When $G$ is simply laced, we can 
give a product expansion for the generating function of the Brylinski 
filtration on the $G$-invariants in the \textit{basic representation} 
$\bfH_0$, the highest-weight module of level $1$ and highest weight $0$.
\begin{theorem}\label{11-12}
For simply laced $G$, the vacuum vector $\omega\in\bfH_0$ gives an isomorphism 
\[\tag{*}
\omega\otimes:\left\{\sym^p(\frg[[z]]dz)^*\right\}^{\frg[z]}
\xrightarrow{\,\sim\,}
\left\{\bfH_0\otimes\sym^p(\frg[[z]]dz)^*\right\}^{\frg[z]}.
\]
Consequently, with $q = z^{-1}$,
\[
\sum_{p\ge 0} t^p\dim_q\mathrm{Gr}_p\,\bfH_0^G = 
	\prod_{k=1}^\ell \prod_{\;n> m_k}(1 - t^{m_k+1}q^n)^{-1}.
\]
\end{theorem}
\begin{proof}
After summing over $p$, the $q$-dimension of the left side in (*) is 
$\prod_{k=1}^\ell\prod_{n>m_k}(1-q^n)^{-1}$ (Theorem \ref{sym}). According 
to \cite{frk}, Theorem ?, this is also the $q$-dimension of $\bfH_0^G$. 
However, the map (*) is an inclusion; hence, using Theorem \ref{11-8}, 
it is an isomorphism, and then it is so in each $p$-degree separately. 
\end{proof}


\section{Line bundle twists} \label{twist}
Let $G$ be simple and simply connected and call $\cO(h)$ the level 
$h$ line bundle on $X$ or $\bfX_\Sigma$. The loop group acts projectively 
on $\cO(h)$, and hence on its Dolbeault cohomologies $H^q(X;\Omega^p(h))$, 
which turn out to be duals of integrable highest-weight modules at level 
$h$, direct products of their $z$-weight spaces. (This follows as in 
Prop.~\ref{62}, except that the cohomologies of $\Gr^n \Omega^p(h)$ 
are now finite sums of duals of irreducible highest-weight modules; 
this suffices for the Mittag-Leffler conditions, as their $z$-graded 
components are finite-dimensional.) For thick flag varieties, we 
obtain instead sums of highest-weight modules (\cite{tel3}, Remark 
8.10). 

The Dolbeault groups of $\cO(h))$ also assemble to a bi-graded module over 
the Dolbeault algebra $H^\bullet(\Omega^\bullet)$. For simply laced $G$ at 
level $1$, our knowledge of the basic invariants (Theorem\ \ref{11-12}) allows 
us to describe the entire structure: $H^\bullet(\Omega^\bullet(h))$ is the 
free module generated from $H^0\left(X;\cO(1)\right)$ under the action of 
the \textit{odd} Dolbeault generators. We prove the theorem for the thick 
loop Grassmannian $\bfX= LG/G[z^{-1}]$; the thin $X$ can be handled as in 
\S\ref{97}. Note that $\bfX=\bfX_{\bA^1}$, with coordinate $q=z^{-1}$ on 
$\bA^1 = \bP^1\setminus\{0\}$.
 
\begin{maintheorem}\label{level1}
For simply laced $G$, $H^\bullet(\bfX,\Omega^\bullet(1))$ is freely 
generated from $\bfH_0= H^0\left(\bfX;\cO(1)\right)$ by the cup-product 
action of the odd generators $\bC[q]dq\subset H^m(\bfX,\Omega^{m+1})$, 
$m= m_1,\ldots,m_\ell$. The multiplication action of the even 
generators of $H^\bullet(\Omega^\bullet)$ is nil.
\end{maintheorem}

\begin{proof}
The centre of $G$ acts trivially on the $H^q(\Omega^p)$; for simply 
laced groups at level $1$, this only allows the basic representation 
$\bfH_0$. The argument now parallels the level zero case, and uses the 
module structure over the latter. By cohomology vanishing (Theorem 
\ref{vanish}), the $E_1^{k,l}$ term replacing \eqref{94} in the sequence 
converging to $H^{k+l}\left(\frM({\bP^1}); \Omega^r(1)\right)$ is now
\[
H^0_{G[[z]]}\left(H^l\left(\bfX; \Omega^{r-k}(1)\right) 
		\otimes\sym^k\frg[[z]]^*\right) \cong 
		H^l\left(\bfX; \Omega^{r-k}(1)\right)^{G[[z]]}
		\otimes H^0_{G[[z]]}\left(\sym^k\frg[[z]]^*\right), 
\]
where we have used the isomorphism of Theorem \ref{11-12}. According to 
\cite{tel4}, Theorem 7.1, the Dolbeault cohomology $H^\bullet\left(
\frM({\bP^1}); \Omega^\bullet(1)\right)$ is isomorphic to $H^{\bullet,
\bullet}(BG;\bC)$, under restriction to the semi-stable sub-stack $BG$ 
of $\frM(\bP^1)$. The argument of \S\ref{proof} now shows that our new 
sequence is freely generated by $\bfH_0$ over the second family of 
level $0$ generators in Proposition \ref{95}.
\end{proof}

\begin{remark} 
This result has an obvious analogue, with parallel proof, for the thick 
flag varieties $\bfX_\Sigma$, when $\Sigma$ has genus 0. Extension to 
higher genus would require us to equate $H^\bullet\left(\frM(\oSig);
\Omega^\bullet(1)\right)$ with the free module spanned by $H^0\left(
\frM(\oSig);\cO(1)\right)$ on half the generators of $H^{p,q}\left(
\frM(\oSig)\right)$. While we believe that equality holds, additional 
arguments seem to be needed.   
\end{remark}

\subsection{Affine Hall-Littlewood functions.} \label{afhall}
For a $G$-representation $V$ with associated vector bundle $\cV$ on 
$\bfX$ (\S\ref{11-7}), the series of characters for the $G$-translation and 
the $z$-scaling actions
\begin{equation}\label{hallit}
P_{h,V}(q,t):= \sum\nolimits_{r,s} (-1)^s (-t)^r \mathrm{ch}
		H^s(\bfX,\Omega^r(h)\otimes \cV)\in R_G[[q,t]]
\end{equation}
are affine analogues of the Hall-Littlewood symmetric functions.\footnote
{The affine Hall-Littlewood functions involve the full flag variety $LG/
\exp(\frB)$ in lieu of the loop Grassmannian, but there is a close 
relation between the two.} We can decompose the $H^q(\bfX;\Omega^p(h))$ 
into the irreducible characters at level $h$, with cofactors $\langle 
P_{h,V} |\bfH\rangle(q,t)\in\bZ[[q,t]]$:
\[
P_{h,V}(q,t) = \sum\nolimits_\bfH \langle P_{h,V} |\bfH\rangle(q,t)
		\cdot \mathrm{ch}(\bfH).
\]
Thus, for simply laced $G$ at level $1$, Theorem \ref{level1} gives for 
the trivial representation $V=\bC$ 
\begin{equation}\label{hall1}
\langle P_{h,\bC} |\bfH_0\rangle(q,t) = \prod_{k=1}^\ell\prod_{n>0}
		(1-t^{m_k+1}q^n).
\end{equation}
Little seems to be known about the cohomology of $\Omega^p(h)\otimes \cV$ 
in general, but the $\langle P_{h,V} |\bfH\rangle(q,t)$ are closely related 
to the Brylinski filtration of \S\ref{11}, as we now illustrate in a simple 
example.

\subsection{Hall-Littlewood cofactors and Brylinski filtration.} For any 
simply connected $G$, the spectral sequence of \S\ref{91} becomes, at 
level $h>0$
\[
E_1^{k,l} = \sum\nolimits_\bfH
	\left\langle H^l\left(\bfX; \Omega^{r-k}(h)\right)|\bfH\right\rangle
	\cdot \left\{\bfH\otimes\sym^k\frg[[z]]^*\right\}^{G[[z]]}
	\Rightarrow H^{k+l}(BG;\Omega^r),
\] 
because $H^\bullet\left(\frM_{\bP^1};\Omega^r(h)\right) = H^\bullet(BG; 
\Omega^r)$. We now form the $(q,t)$-characteristic by multiplying the 
left side by $(-1)^{k+l}(-t)^r$ and summing over $k,l,r$. Theorem \ref{11-8} 
(with the substitution $t\mapsto tq^{-1}$, to undo the shift introduced 
in \S\ref{11-1}) gives the near-orthogonality relation 
\begin{equation}\label{ortho}
\sum\nolimits_\bfH \langle P_{h,\bC}|\bfH\rangle(q,t)\cdot
	\sum\nolimits_p (tq^{-1})^p\dim_q\mathrm{Gr}_p\,\bfH^G =
	\prod\nolimits_{k=1}^\ell (1-t^{m_k+1})^{-1};
\end{equation}
the right-hand side is $\sum_{r,s}(-1)^s(-t)^r h^s(BG;\Omega^r) = \sum_r 
t^r h^{2r}(BG)$. Implications of \eqref{ortho} will be explored in 
future work; instead, we conclude with a combinatorial application.

\subsection{Lattice hyper-geometric sums.}\label{bailey}
There is a Kac formula for $P_{1,\bC}$, established as in \S\ref{euler} 
(but now with the \textit{increasing} filtration on $\Omega^p$, as we work 
on the thick Grassmannian $\bfX$):
\begin{equation}\label{waffaction}
\sum_{w\in\waff} w\left[\prod_{n>0;\,\alpha}
		\frac{1-tq^n\me^\alpha}{1-q^n\me^\alpha}
		\cdot\prod_{\alpha>0}(1-\me^\alpha)^{-1}\right]
		\cdot\prod_{n>0}\left(\frac{1-tq^n}{1-q^n}\right)^\ell. 
\end{equation} 
At level $1$, a lattice element $\gamma\in\waff$ sends $q^n\me^\lambda$ 
to $q^{n+\gamma^2/2+\langle \lambda |\gamma\rangle}\me^{\lambda+\gamma}$, 
in which the basic inner product is used to convert $\gamma$ to a weight. 
The manipulation in \S\ref{1psi1} converts \eqref{waffaction} into
\begin{equation}
\sum_{\gamma} q^{\gamma^2/2}\me^{\gamma}\cdot
	\prod_{n>0;\,\alpha}
		\frac{1-tq^{n+\langle\alpha|\gamma\rangle} \me^\alpha}
			{1-q^{n+\langle\alpha|\gamma\rangle}\me^\alpha} 
	\cdot\prod_{n>0}\left(\frac{1-tq^n}{1-q^n}\right)^\ell.
\end{equation}
For simply laced $G$, another expression is provided by \eqref{hall1} and 
any of the character formulae for $\bfH_0$; thus, the basic bosonic 
realisation gives   
\begin{equation}\label{12-10}
P_{1,\bC} = \prod_{k=1}^\ell \prod_{n>0} \frac{1-t^{m_k+1}q^n}{1-q^n}
	\cdot\sum_{\gamma} q^{\gamma^2/2}\me^\gamma
\end{equation}
where we sum over the co-root lattice (which is also the root lattice). 
Equating the last two expressions gives the identity
\begin{equation}\label{12-11}
\sum_{\gamma} q^{\gamma^2/2}\me^{\gamma}\cdot
	\prod_{n>0;\,\alpha}
		\frac{1-tq^{n+\langle\alpha|\gamma\rangle} \me^\alpha}
			{1-q^{n+\langle\alpha|\gamma\rangle}\me^\alpha} 
	= \prod_{k=1}^\ell\prod_{n>0} \frac{1-t^{m_k+1}q^n}{1-tq^n}
		\cdot\sum_{\gamma} q^{\gamma^2/2}\me^\gamma.
\end{equation}
With $G=\SL_2$, replacing $q$ by $\sqrt q$ leads to 
\[
\sum_{m\in\bZ} q^{m^2/2}u^{2m}\cdot\prod_{n>0}
		\frac{(1-tq^{n/2+m} u^2)(1-tq^{n/2-m} u^{-2})}
			{(1-q^{n/2+m}u^2)(1-q^{n/2-m}u^{-2})} 
	= \prod_{n>0} \frac{1-t^2q^{n/2}}{1-tq^{n/2}}
		\cdot\sum_{m\in\bZ} q^{m^2/2}u^{2m},
\]
which, using the notation $(a)_\infty=\prod_{n\ge 0}(1-aq^n)$, 
$(a)_n = (a)_\infty/(aq^n)_\infty$, $(a_1,\dots,a_k)_n = \prod_i(a_i)_n$, 
becomes the hyper-geometric summation formula
\begin{multline*}
\sum_{m\in\bZ} q^{m^2/2}u^{2m}\cdot\frac{(\sqrt qu^2)_{m}
	(\sqrt qu^{-2})_{-m}(qu^2)_{m}(qu^{-2})_{-m}}
	{(\sqrt qtu^2)_{m}(\sqrt qtu^{-2})_{-m}(qtu^2)_{m}(qtu^{-2})_{-m}}\\ 
	= \frac{(\sqrt qu^2,\sqrt qu^{-2},qu^2,qu^{-2},\sqrt q t^2,	
		q t^2)_\infty}
		{(\sqrt qtu^2,\sqrt qtu^{-2},qtu^2,qtu^{-2},\sqrt qt,qt)_\infty}
		\cdot\sum_{m\in\bZ} q^{m^2/2}u^{2m};
\end{multline*}
most factors in the numerator on the left cancel out, and the series can 
then be summed by specialising Bailey's $_4\psi_4$ summation formula 
(see \cite{gasp}, Ch.\ 5).

\vskip1cm
\part*{Appendix}
\addtocontents{toc}{\small\bfseries Appendix}

\appendix
\section{Proof of Lemma \ref{313}}\label{app}
It is clear that both sides in \eqref{313} annihilate the constant line
in $\Lambda \otimes \sym$, and it is also easy to see that they agree on 
the symmetric part $1\otimes\sym$, where $D$, $K$, $\ad$ and $\ad^*$ 
all vanish. So we must only check equality on the linear $\psi$ terms, 
and on the quadratic $\psi\wedge\psi$ and $\sigma \psi$ terms. 

\subsection{The linear $\psi$ terms.} \label{A1}
Fix $b\in A$, $n>0$. We compute: 
\[
\bar\partial \psi^b(-n) = \frac{1}{2}\sum_
	{ \ofrac{a\in A}{ {0<m<n}}}
		 {\psi^a(-m)\wedge\psi^{[a,b]}(m-n)},
\] 
\begin{align*}
\bar\partial^*\bar\partial \psi^b(-n) 
&= \frac{1}{4} \sum_{\ofrac{a\in A} { {0<m<n}}}
			{\frac{n}{m(n-m)}\psi^{[a,[a,b]]}(-n)} - 
	\frac{1}{4}\sum_{\ofrac{a\in A}{0<m<n}} 
			{\frac{n}{m(n-m)}\psi^{[[a,b],a]}(-n)} \\ 
&= \frac{1}{2}\sum_{ \ofrac{a\in A} { {0<m<n}}} 
		\frac{n}{m(n-m)}\psi^{[a,[a,b]]}(-n) = 
			\frac{1}{2}\sum_{0<m<n} {\frac{n}{m(n-m)}\psi^b(-n)}\\
&= \frac{1}{2}\sum_{0<m<n} {\left(\frac{1}{m} + \frac{1}{n-m}
		\right)\cdot \psi^b(-n)} = \sum_{0<m<n} {1/m\cdot \psi^b(-n)}.
\end{align*}
Further, $\bar\partial^*\psi^b(-n)=0$, so $\overline\square\psi^b(-n)$ is 
as just computed. Next, 

\begin{align*}
\square\psi^b(-n) &=\sum_{\ofrac{a\in A}{0<m<n}} 
		\frac{1}{m}\ad_a(-m)\ad_a(-m)^*\psi^b(-n)
	=\sum_{ \ofrac{a\in A} { {0<m<n}}}
\frac{1}{m}\cdot \frac{n-m}{n}\cdot \psi^{[a,[a,b]]}(-n) \\
& =\sum_{0<m<n}\left(\frac{1}{m} - \frac{1}{n}\right)\cdot 
		\psi^b(-n) = \sum_{0<m<n} \frac{1}{m}\cdot \psi^b(-n)
				- \psi^b(-n) + \psi^b(-n)/n \\					
D\psi^b(-n) &= \sum\limits_{\ofrac{a\in A}{0<m<n}} 
		\psi^{[a,[a,b]]}(-n)/n = \psi^b(-n),\\
K\psi^b(-n) &= \frac{1}{n} \sum_{a\in A}{\ad_{[a,b]}(0)\psi^a(-n)} 
			= -\frac{1}{n}\sum_{a\in A} \psi^{[a,[a,b]]}(-n) 
			= -\psi^b(-n)/n,
\end{align*} 
and the last three terms sum up to $\overline\square\psi^b(-n)$, as claimed. 
\qed 
 
\subsection{The quadratic $\psi\wedge\psi$ terms.}\label{A2}
Fix $b,c\in A$ and $n,p>0$. For each second-order differential operator $P$ 
involved, we focus on the  \textit{cross-term} $P\left(\psi^b(-n)\wedge 
\psi^c(-p)\right) - P\psi^b(-n)\wedge \psi^c(-p) - \psi^b(-n)\wedge 
P\psi^c(-p)$; equality of cross-terms and the identity \eqref{313} on the 
linear factors imply the identity for quadratic terms. 

\begin{lemma}\label{A3}
The cross-term in $\overline\square\left(\psi^b(-n)\wedge\psi^c(-p)\right)$ 
is the following three-term sum: 
\begin{align*}
&\sum_{\ofrac{a\in A}{0<m<n}} 
	\left(\frac{1}{m} - \frac{1}{n}\right)\cdot 
			\psi^{[a,b]}(m-n)\wedge\psi^{[a,c]}(-m-p)\\
+& \sum_{\ofrac{a\in A} {0<m<p}}
	\left(\frac{1}{m} - \frac{1}{p}\right)\cdot 
		\psi^{[a,b]}(-m-n)\wedge\psi^{[a,c]}(m-p)\\
-& \sum_{a\in A} \left(\frac{1}{n} + \frac{1}{p}\right)\cdot
			\psi^{[a,b]}(-n)\wedge\psi^{[a,c]}(-p).
\end{align*}
\end{lemma}
\begin{proof} We rewrite the sum by adding and subtracting terms: 
\begin{align}
&\sum_{\ofrac{a\in A} {0<m<n}} 
	\left(\frac{1}{m} + \frac{1}{p}\right) \cdot 
		\psi^{[a,b]}(m-n)\wedge \psi^{[a,c]}(-m-p) \nonumber\\
+& \sum_{ \ofrac{a\in A} {0<m<p}} 
\left(\frac{1}{m} + \frac{1}{n}\right)\cdot 
		\psi^{[a,b]}(-m-n)\wedge \psi^{[a,c]}(m-p) \label{A4} \\ 
-& \sum_{ \ofrac{a\in A} {0<m<n+p}}
		\left(\frac{1}{n} + \frac{1}{p}\right)\cdot 
			\psi^{[a,b]}(-m)\wedge\psi^{[a,c]}(m-n-p).\nonumber
\end{align}
We now track, in turn, the source of each of the three terms in \eqref{A4}. 
We have 
\begin{multline}\label{A5}
\bar\partial\left(\psi^b(-n)\wedge\psi^c(-p)\right)
		= \frac{1}{2}\sum_{\ofrac{a\in A} {0<m<n}}
			\psi^a(-m)\wedge \psi^{[a,b]}(m-n)\wedge \psi^c(-p)\\
		+ \frac{1}{2}\sum_{\ofrac{a\in A} {0<m<p}} 
			\psi^a(-m)\wedge\psi^b(-n)\wedge\psi^{[a,c]}(m-p),
\end{multline}
and applying $\bar\partial^*$ to the first sum gives the following, 
after collecting the terms where $\psi^c(-p)$ survives intact into the
first summand: 
\begin{align}\label{A6}
&\bar\partial^*\bar\partial \psi^b(-n)\wedge \psi^c(-p)\nonumber\\
&+ \frac{1}{4}\sum_{\ofrac{a\in A} {0<m<n}}
	\frac{1}{m}\psi^{[a,b]}(m-n)\wedge \ad_a(m)^* \psi^c(-p)\nonumber\\
&- \frac{1}{4}\sum_{\ofrac{a\in A} {0<m<n}}
	\frac{1}{n-m} \psi^a(-m)\wedge \ad_{[a,b]}(n-m)^* \psi^c(-p)\\
&+ \frac{1}{4}\sum_{\ofrac{a\in A} {0<m<n}}
	\frac{1}{p}\ad_c(p)^*\left(\psi^a(-m)\wedge\psi^{[a,b]}(m-n)\right).
	\nonumber
\end{align}
The first term is $\overline\square\psi^b(-n)\wedge \psi^c(-p)$. The two
middle-line terms agree, after substituting $m\leftrightarrow (n-m)$, and 
sum to 
\begin{equation}\label{A7}
\frac{1}{2}\sum_{\ofrac{a\in A} {0<m<n}} 
	\frac{p+m}{mp}\psi^{[a,b]}(m-n)\wedge \psi^{[a,c]}(-m-p).
\end{equation} 
Amusingly, the third line takes the same value \eqref{A7}; so the sum in 
\eqref{A6} equals 
\begin{equation}\label{A8}		
\overline\square\psi^b(-n)\wedge \psi^c(-p)
	+\sum_{\ofrac{a\in A} {0<m<n}}
		\frac{p+m}{mp}\psi^{[a,b]}(-m-p)\wedge \psi^{[a,c]}(m-n),
\end{equation} 
and so the cross-term in \eqref{A8} accounts for the first term in 
\eqref{A4}. Substituting $b\leftrightarrow c$, $n\leftrightarrow p$ shows 
that the $\bar\partial^*$-image of the second term in \eqref{A5} is 
\begin{equation}\label{A9}		
\psi^b(-n)\wedge\overline\square\psi^c(-p) + 
	\sum_{\ofrac{a\in A} {0<m<n}}
		\frac{n+m}{mn}\psi^{[a,b]}(-m-n)\wedge \psi^{[a,c]}(m-p),
\end{equation} 
whose cross-term is the second term in \eqref{A4}. 

Finally, $\bar\partial^*\left(\psi^b(-n)\wedge \psi^c(-p)\right) = \frac{n+p}
{np}\cdot \psi^{[b,c]}(-p-n)$, whence by applying $\bar\partial$ we get
\begin{multline*}
\bar\partial\bar\partial^*\left(\psi^b(-n)\wedge \psi^c(-p)\right)
	= \frac{1}{2}\frac{n+p}{np} 
		\sum_{\ofrac{a\in A} {0<m<n+p}}
			\psi^a(-m)\wedge \psi^{[a,[b,c]]}(m-p-n) \\ 
	= \frac{n+p}{2np} 
		\sum\limits_{\ofrac{a\in A}{0<m<n+p}}
		 	\left(\psi^a(-m) \wedge \psi^{[[a,b],c]}(m-p-n) +
			\psi^a(-m)\wedge \psi^{[b,[a,c]]}(m-p-n)\right)\\
	= -\frac{n+p}{np}
		\sum_{\ofrac{a\in A} {0<m<n+p}}
			\psi^{[a,b]}(-m)\wedge \psi^{[a,c]}(m-p-n),
\end{multline*}
which is the third term in \eqref{A4}. The proposition is proved. 
\end{proof}

Now $D\left(\psi^b(-n)\wedge\psi^c(-p)\right) = D\psi^b(-n)\wedge\psi^c(-p)
+ \psi^b(-n)\wedge D\psi^c(-p)$, with no cross-term, while
\begin{multline}\label{A10}
K\left(\psi^b(-n)\wedge \psi^c(-p)\right)
	= K\psi^b(-n)\wedge \psi^c(-p) + \psi^b(-n)\wedge K\psi^c(-p)\\
- \frac{1}{n}\psi^{[a,b]}(-n)\wedge \psi^{[a,c]}(-p) - \frac{1}{p}
	\psi^{[a,b]}(-n)\wedge \psi^{[a,c]}(-p), 
\end{multline}
\begin{align}\label{A11}
\square\left(\psi^b(-n)\wedge\psi^c(-p)\right) 
	&= \psi^b(-n)\wedge \psi^c(-p) + \psi^b(-n)\wedge\psi^c(-p)\nonumber\\
	&+ \sum_{\ofrac{a\in A} {0<m<n}}
	\frac{1}{m}\frac{n-m}{n}\psi^{[a,b]}(m-n)\wedge\psi^{[a,c]}(-m-p)\\
	&+ \sum_{\ofrac{a\in A} {0<m<p}}
	\frac{1}{m}\frac{p-m}{p}\psi^{[a,b]}(-m-n)\wedge \psi^{[a,c]}(m-p),
	\nonumber
\end{align}
and the cross-terms in \eqref{A10} and \eqref{A11} add up to the expression 
in \eqref{A3}.

\subsection{The $\sigma\psi$ terms.}
As before, fix $b,c\in A$ and $n,p > 0$. Then, 
\begin{align*}
\bar\partial\left(\sigma^b(-n)\psi^c(-p)\right) 
	&= \sum\limits_{\ofrac{a\in A} {0<m\le n}}
		\sigma^{[a,b]}(m-n)\cdot \psi^a(-m)\wedge \psi^c(-p)\\ 
	&+ \frac{1}{2}\sum_{\ofrac{a\in A} {0<m<p}}
		\sigma^b(-n)\psi^a(-m)\wedge \psi^{[a,c]}(m-p),
\end{align*} 
and applying $\bar\partial^*$ yields the following sum: 
\begin{align*}
&\overline\square \sigma^b(-n)\cdot \psi^c(-p) - \frac{1}{p}
		\sum_{\ofrac{a\in A} {0<m\le n}} 
			\sigma^{[a,c]}(m-p-n)\psi^{[a,b]}(-m)\\ 
&+ \frac{1}{2} \sum_{\ofrac{a\in A}{0<m<p}}
			\frac{1}{m}\sigma^{[a,b]}(-m-n)\psi^{[a,c]}(m-p) 
			- \frac{1}{2} \sum_{\ofrac{a\in A} {0<m<p}}
			\frac{1}{m-p}\sigma^{[[a,c],b]}(m-p-n)\psi^a(-m)\\ 
&+ \frac{1}{2}\sum_{\ofrac{a\in A}{0<m\le n}}
		\frac{m+p}{mp}\sigma^{[a,b]}(m-n)\psi^{[a,c]}(-m-p) +
		\frac{1}{2}\sum_{\ofrac{a\in A} {0<m\le n}}
			\frac{m+p}{mp}\sigma^{[a,b]}(m-n)\psi^{[a,c]}(-m-p)\\
&+ \sigma^b(-n)\overline\square\psi^c(-n).
\end{align*}
The first two lines come from the $R^*$-terms in $\bar\partial^*$, the 
last two lines from the $\ad^*$-terms. The two terms in each of the middle 
rows are equal, so the cross-term can be rewritten as follows: 
\begin{equation}\begin{split}\label{A13}
-\frac{1}{p}\sum_{\ofrac{a\in A}{0<m\le n}}
	\sigma^{[a,c]}(m-p-n)\cdot \psi^{[a,b]}(-m) 
&+ \sum_{\ofrac{a\in A} {0<m<p}}
		\frac{1}{m}\sigma^{[a,b]}(-m-n)\cdot \psi^{[a,c]}(m-p) \\
+ \sum_{\ofrac{a\in A}{0<m\le n}} 
	\left(\frac{1}{m} + \frac{1}{p}\right)
		&\sigma^{[a,b]}({m-n} )\cdot \psi^{[a,c]}(-m-p)
\end{split}\end{equation}
Now, $\bar\partial^*\left(\sigma^b(-n)\cdot\psi^c(-p)\right) = 
\sigma ^{[c,b]}(-n-p)/p$, whence
\begin{multline}\label{A14}
\bar\partial\bar\partial^*\left(\sigma^b(-n)\cdot \psi^c(-p)\right)
	= \frac{1}{p}\sum_{\ofrac{a\in A} {0<m\le n+p}}
		\sigma^{[a,[c,b]]}(m-n-p)\cdot\psi^a(-m) \\ 
	= \frac{1}{p}\sum_{\ofrac{a\in A}{0<m<n+p}} 
		\sigma ^{[a,c]}(m-n-p)\cdot \psi^{[a,b]}(-m) 
	- \frac{1}{p}\sum_{\ofrac{a\in A} {0<m<n+p}} 
		\sigma ^{[a,b]}(m-n-p)\cdot \psi^{[a,c]}(-m).
\end{multline}
Summing \eqref{A13} and \eqref{A14} gives
\begin{multline}\label{A15}
\frac{1}{p}\sum_{\ofrac{a\in A}{0<m\le p}} 
	{\sigma^{[a,c]}(m-p)\cdot \psi^{[a,b]}(-m-n)} +
		\sum\limits_{\ofrac{a\in A}{0<m<p}} 
		\left(\frac{1}{ m} - \frac{1}{p}\right)\sigma^{[a,b]}(-m-n)\cdot 
				\psi^{[a,c]}(m-p)\\
+ \sum_{\ofrac{a\in A} {0<m<n}} 
	\frac{1}{m}\sigma^{[a,b]}(m-n)\cdot \psi^{[a,c]}(-m-p) -
	\frac{1}{p}\sum_{a\in A} \sigma^{[a,b]}(-n)\cdot \psi^{[a,c]}(-p);
\end{multline}
here, the first term is the sum of the first terms in \eqref{A13} and 
\eqref{A14}, the second and third incorporate the second and third terms in 
\eqref{A13} and the $0<m<p$, resp.\ the $p<m<p+n$ portions of the second term 
in \eqref{A14}, and the final term is the $m=p$ contribution of the same. 

Moving on to the right-hand side of (\ref{313}), the cross-term in 
$\square\left(\sigma^b(-n)\psi^c(-p)\right)$ is 
\begin{equation}\label{A16}
\sum_{\ofrac{a\in A} {0<m<n}}
	\frac{1}{m} \sigma^{[a,b]}(m-n)\cdot \psi^{[a,c]}(-m-p)
+\sum_{\ofrac{a\in A} {0<m<p}}
	\frac{p-m}{mp}\sigma^{[a,b]}(-m-n)\cdot \psi^{[a,c]}(m-p),
\end{equation} 
the two terms coming from the $\ad\cdot R^*$ and $R\cdot\ad^*$ cross-terms, 
respectively. Further,
\begin{equation}\label{A17}
D\left(\sigma^b(-n)\cdot \psi^c(-p)\right) = 
		\sigma^b(-n)\cdot D\psi^c(-p) + \frac{1}{p}
			\sum_{\ofrac{a\in A}{0<m\le p}} 
				\sigma^{[a,c]}(m-p)\cdot\psi^{[a,b]}(-m-n),		
\end{equation}
\begin{equation}\label{A18}
K\left(\sigma^b(-n)\cdot \psi^c(-p)\right) = 
	\sigma^b(-n)\cdot K\psi^c(-p) - \frac{1}{p} \sum_{a\in A} 
		\sigma^{[a,b]}(-n)\cdot \psi^{[a,c]}(-p).
\end{equation} 
It is now clear that the cross-terms in \eqref{A16}--\eqref{A18} sum 
to \eqref{A15}.  \qed 

\vskip2cm

\vspace{1.5cm}
\small{\noindent
\textsc{Susanna Fishel:} Mental Images, Fasanenstra\ss e 81, 10623 Berlin,
Germany \texttt{susanna@mental.com}\\
\textsc{Ian Grojnowski, Constantin Teleman:} DPMMS, Wilberforce Road, 
Cambridge CB3 0WB, UK \\
\texttt{groj@dpmms.cam.ac.uk}, \texttt{teleman@dpmms.cam.ac.uk}}

\end{document}